\documentclass[12pt]{amsart}
\usepackage{amssymb,amsmath,verbatim,fullpage}
\makeatletter
\def\@cite#1#2{{\m@th\upshape\bfseries%
[{#1\if@tempswa{\m@th\upshape\mdseries, #2}\fi}]}}
\makeatother
\theoremstyle{plain}
\newtheorem{thm}{Theorem}[section]
\newtheorem{cor}[thm]{Corollary}
\newtheorem{prop}[thm]{Proposition}
\newtheorem{lem}[thm]{Lemma}

\theoremstyle{definition}
\newtheorem{rem}[thm]{Remark}
\newtheorem{defn}[thm]{Definition}

\newtheorem{eg}[thm]{Example}

\newtheorem{ques}[thm]{Question}
\newcommand{\Prf}{\noindent\textbf{Proof.\ }}

\renewcommand{\qed}{\hfill \vrule height5pt width5pt depth1pt \vspace{+2.00ex}}
\newcommand{\upqed}{\vspace{-2.5\baselineskip}\newline\hbox{}\qed}
\newcommand{\bsl}{\setminus}

\newcommand{\ca}{\mathrm{C}^*}

\newcommand{\dlim}{\displaystyle\lim\limits}
\newcommand{\dsum}{\displaystyle\sum\limits}
\newcommand{\tsum}{\textstyle\sum\limits}
\newcommand{\Fn}{\bF_n^+}
\newcommand{\Hn}{\ltwo(\bF_n^+)}
\newcommand{\lip}{\langle}
\newcommand{\rip}{\rangle}
\newcommand{\ip}[1]{\lip #1 \rip}

\newcommand{\mt}{\varnothing}
\newcommand{\norm}[1]{\left\| #1 \right\|}
\newcommand{\ol}{\overline}

\newcommand{\upminus}{\raise.9ex\hbox{-\!}}

\newcommand{\sot}{\textsc{sot}}
\newcommand{\wot}{\textsc{wot}}
\newcommand{\wotclos}[1]{\ol{#1}^{\textsc{wot}}}

\newcommand{\wsclos}[1]{\ol{#1}^{\text{w-}*}}

\DeclareMathOperator*{\sotlim}{\textsc{sot}--lim}

\newcommand{\sotsum}{\textsc{sot--}\!\!\sum}

\newcommand{\vac}{\V_{ac}}

\newcommand{\bA}{{\mathbb{A}}}

\newcommand{\bC}{{\mathbb{C}}}
\newcommand{\bD}{{\mathbb{D}}}

\newcommand{\bF}{{\mathbb{F}}}

\newcommand{\bN}{{\mathbb{N}}}

\newcommand{\bT}{{\mathbb{T}}}

\newcommand{\bZ}{{\mathbb{Z}}}
  \newcommand{\A}{{\mathcal{A}}}
  \newcommand{\B}{{\mathcal{B}}}

  \newcommand{\E}{{\mathcal{E}}}
  
  \newcommand{\G}{{\mathcal{G}}}
\renewcommand{\H}{{\mathcal{H}}}

  \newcommand{\M}{{\mathcal{M}}}
  \newcommand{\N}{{\mathcal{N}}}
\renewcommand{\O}{{\mathcal{O}}}

\renewcommand{\S}{{\mathcal{S}}}
  
  \newcommand{\U}{{\mathcal{U}}}
  \newcommand{\V}{{\mathcal{V}}}
  \newcommand{\W}{{\mathcal{W}}}
  \newcommand{\X}{{\mathcal{X}}}
  \newcommand{\Y}{{\mathcal{Y}}}

\newcommand{\ep}{\varepsilon}
\newcommand{\eps}{\varepsilon}
\renewcommand{\phi}{\varphi}

\newcommand{\upchi}{{\raise.35ex\hbox{$\chi$}}}

\newcommand{\fA}{{\mathfrak{A}}}

\newcommand{\fJ}{{\mathfrak{J}}}
\newcommand{\fK}{{\mathfrak{K}}}
\newcommand{\fL}{{\mathfrak{L}}}
\newcommand{\fM}{{\mathfrak{M}}}

\newcommand{\fR}{{\mathfrak{R}}}
\newcommand{\fS}{{\mathfrak{S}}}
\newcommand{\fs}{{\mathfrak{s}}}
\newcommand{\fT}{{\mathfrak{T}}}
\newcommand{\fU}{{\mathfrak{U}}}

\newcommand{\fW}{{\mathfrak{W}}}
\newcommand{\fX}{{\mathfrak{X}}}

\newcommand{\AD}{\mathrm{A}(\mathbb{D})}

\newcommand{\ltwo}{\ell^2}

\newcommand{\AND}{\text{ and }}

\newcommand{\qand}{\quad\text{and}\quad}

\newcommand{\qforal}{\quad\text{for all}\quad}
\newcommand{\qif}{\quad\text{if}\quad}

\newcommand{\Alg}{\operatorname{Alg}}

\newcommand{\diag}{\operatorname{diag}}
\newcommand{\Dim}{\operatorname{dim}}

\newcommand{\dist}{\operatorname{dist}}

\newcommand{\im}{\operatorname{Im}}

\newcommand{\Lat}{\operatorname{Lat}}
\newcommand{\ran}{\operatorname{Ran}}

\newcommand{\re}{\operatorname{Re}}
\newcommand{\Rep}{\operatorname{Rep}}

\newcommand{\spn}{\operatorname{span}}

\begin{document}

%%%%%%%%%%%%%%%%%%%%%%%%%%%%%%%%%%%%%%%%%%%%%%%%%%%%%%%%%%%%
\title[Free Semigroup Algebras]%
{Absolutely Continuous Representations\\
 and a Kaplansky Density Theorem\\
for Free Semigroup Algebras}
\thanks{Printed on \today.}
\author[K.R.Davidson]{Kenneth R. Davidson}
\thanks{First author partially supported by an NSERC grant}
\address{Pure Math.\ Dept.\\U. Waterloo\\Waterloo, ON\;
N2L--3G1\\CANADA}
\email{krdavids@uwaterloo.ca}
\author[J.Li]{Jiankui~Li}
%\thanks{Second author's research partially supported by ?????}
\address{Math.\ Dept.\\E.C.U.S.T.\\Shanghai 200237 \\ P.R. CHINA }
\email{jiankuili@yahoo.com}
\author[D.R.Pitts]{David~R.~Pitts}
\address{Math.\ Dept.\\University of Nebraska\\
Lincoln, NE 68588\\USA}
\email{dpitts@math.unl.edu}
%\thanks{Third author partially supported by an NSF grant}
%
\subjclass{47L80}
\date{}

\begin{abstract}
 We introduce notions of absolutely continuous
functionals and representations on the non-com\-mu\-ta\-tive disk
algebra $\fA_n$.  Absolutely continuous functionals are used to help
identify the type L part of the free semigroup algebra associated to a
$*$-extendible representation $\sigma$.  A $*$-extendible
representation of $\fA_n$ is \textit{regular} if the absolutely continuous
part coincides with the type L part.  All known examples are regular.
Absolutely continuous functionals are intimately related to maps which
intertwine a given $*$-extendible representation with the left regular
representation.  A simple application of these ideas extends
reflexivity and hyper-reflexivity results.  Moreover the use of
absolute continuity is a crucial device for establishing a density
theorem which states that the unit ball of $\sigma(\fA_n)$ is weak-$*$
dense in the unit ball of the associated free semigroup algebra if and
only if $\sigma$ is regular.  We provide some explicit constructions
related to the density theorem for specific representations.  A notion
of singular functionals is also defined, and every functional
decomposes in a canonical way into the sum of its absolutely
continuous and singular parts.
\end{abstract}
\maketitle
%%%%%%%%%%%%%%%%%%%%%%%%%%%%%%%%%%%%%%%%%%%%%%%%%%%%%%%%%%%%

Free semigroup algebras were introduced in \cite{DP1} as a method
for analyzing the fine structure of $n$-tuples of isometries with
commuting ranges.  The C*-algebra generated by such an $n$-tuple is
either the Cuntz algebra $\O_n$ or the Cuntz-Toeplitz algebra $\E_n$.  As
such, the free semigroup algebras can be used to reveal the fine spatial
structure of representations of these algebras much in the same way as
the von Neumann algebra generated by a unitary operator encodes the
measure class and multiplicity which cannot be detected in the
C*-algebra it generates.

This viewpoint yields critical information in the work of Bratteli and
Jorgensen \cite{BJ1,BJ3,Jo,JK} who use certain representations of $\O_n$ to
construct and analyze wavelet bases.

{From} another point of view, free semigroup algebras can be used to
study arbitrary (row contractive) $n$-tuples of operators.  Frahzo
\cite{Fr1,Fr2}, Bunce \cite{Bun} and Popescu \cite{Pop_diln} show that
every (row) contractive $n$-tuple of operators has a unique minimal
dilation to an $n$-tuple of isometries which is a row contraction,
meaning that the ranges are pairwise orthogonal. Thus every row
contraction determines a free semigroup algebra.  Popescu
\cite{Pop_vN} establishes the $n$-variable von Neumann inequality
which follows immediately from the dilation theorem.  Popescu has
pursued a program of establishing the analogues of the
Sz.~Nagy--Foia\c{s} program in the $n$-variable setting
\cite{Pop_char,Pop_fact,Pop_func}; the latter two papers deal with
the free semigroup algebras from this point of view.  Free semigroup
algebras play the same role for non-commuting operator theory as the
weakly closed unital algebra determined by the isometric dilation of a
contraction plays for a single operator.

In \cite{DKS}, the first author, Kribs and Shpigel use dilation theory
to classify the free semigroup algebras which are obtained as the
minimal isometric dilation of contractive $n$-tuples of operators on
finite dimensional spaces.  Such free semigroup algebras are called
\textit{finitely correlated}, because from the wavelet perspective,
these algebras correspond to the finitely correlated representations
of $\E_n$ or $\O_n$ introduced and studied by Bratteli and Jorgensen.
It is interesting that this class of representations of $\O_n$ are
understood in terms of an $n$-tuple of matrices, a reversal of the the
single variable approach of analyzing arbitrary operators using the
isometric dilation.  Out of the analysis of finitely correlated free
semigroup algebras emerged a structural result that appeared to rely
on the special nature of the representation.  However in \cite{DKP},
two of the current authors and Katsoulis were able to expose a rather
precise and beautiful structure for arbitrary free semigroup algebras.
This structure plays a key role in this paper.

The prototype for free semigroup algebras is the algebra $\fL_n$
determined by the left regular representation of the free semigroup
$\Fn$ on Fock space.  This representation arises naturally in the
formulation of quantum mechanics.  We have named it the 
\textit{non-commutative analytic Toeplitz algebra} because of the
striking analytic properties that it has
\cite{AP_fact,DP1,DP2,DP3,AP_pois}. In particular, the vacuum vector
(and many other vectors in this representation) has the property that
its image under all words in the
$n$ isometries forms an orthonormal set.  We call such vectors
\textit{wandering vectors}.  Such vectors play a crucial role in these
representations, and a deeper understanding of when they occur is one of
the main open questions in this theory.

The norm-closed algebra $\fA_n$ generated by $n$ isometries with
orthogonal ranges is even more rigid than the C*-algebra.
Indeed, it sits inside the C*-algebra $\E_n$, but the quotient onto
$\O_n$ is completely isometric on this subalgebra.  As $\O_n$ is
simple, it is evidently the C*-envelope of $\fA_n$.
This algebra has been dubbed the \textit{non-commutative disk algebra}
by Popescu.  It plays the same role in this theory as the disk algebra
plays in the study of a single isometry.

In this paper, we explore in greater depth the existence of wandering
vectors.  The major new device is the notion of an \textit{absolutely
continuous} linear functional on $\fA_n$.  In the one variable case, a
functional on $\AD$ is given by integration against a representing
measure supported on the Shilov boundary $\bT$.
Absolute continuity is described in terms of Lebesgue measure.
In our setting, we do not have a boundary, and we have instead defined
absolute continuity in terms of its relationship to the left regular
representation.

A related notion that plays a key role are intertwining maps from the
left regular representation to an arbitrary free semigroup algebra.
The key observation is that the range of such maps span the vectors
which determine absolutely continuous functionals, and they serve to
identify the type L part of the representation (see below).  These
results will be used to clarify precisely when a free semigroup is
reflexive.  For type L representations, we establish hyper-reflexivity
whenever there are wandering vectors---the reflexive case.  Basically
the only obstruction to hyper-reflexivity is the possibility that
there may be a free semigroup algebra which is type L (isomorphic to
$\fL_n$) but has no wandering vectors, and hence will be reductive
(all invariant subspaces have invariant ortho-complements).

The ultimate goal of this paper is to obtain an analogue of the
Kaplansky density theorem.  This basic and well-known result states
that given any C*-algebra and any $*$-representation, the image of the
unit ball is \wot-dense in the unit ball of the \wot-closure.  In the
nonselfadjoint setting, such a result is not generally true.
However, in the context of completely isometric representations of
$\fA_n$, we have a rather rigid structure, and we shall show that in
fact such a Kaplansky type theorem does hold. Let $\sigma$ be a
$*$-extendible representation of $\fA_n$, that is, $\sigma$ is the
restriction of a $*$-representation of $\O_n$ or $\E_n$ to $\fA$.  We
call it regular if the type L part coincides with the absolutely
continuous part.  It is precisely this case in which a density theorem
holds, and the unit ball of $\sigma(\fA_n)$ is weak-$*$ dense in the
unit ball of the free semigroup algebra.  In particular, we shall see
that this holds in the presence of a wandering vector.  In fact, the
only possible obstruction to a Kaplansky density result for all
representations of $\fA_n$ is the existence of a representation where
the free semigroup algebra is a von Neumann algebra and is also
absolutely continuous.  No such representation is known to exist.

%%%%%%%%%%%%%%%%%%%%%%%%%%%%%%%%%%%%%%%%%
\section{Preliminaries}\label{S:prelim}

In this section, we will remind the reader of some of the more
technical aspects which we need, and will establish some
notation for what follows.

A typical $n$-tuple of isometries acting on a Hilbert space $\H$ and
having pairwise orthogonal ranges will be denoted by $S_1,\dots,S_n$. 
This may be recognized algebraically by the relations $S_j^*S_j = I \ge
\sum_{i=1}^n S_iS_i^*$, $1\leq j\leq n.$ 
The C*-algebra that they generate is the Cuntz algebra $\O_n$ when
$\sum_{i=1}^n S_iS_i^* = I$ and the Cuntz-Toeplitz algebra $\E_n$ when
$\sum_{i=1}^n S_iS_i^* < I$.  The norm-closed unital subalgebra
generated by $S_1,\dots, S_n$ (but not their adjoints) is completely
isometrically isomorphic to Popescu's non-commutative disk algebra $\fA_n$.
The ideal of $\E_n$ generated by $I - \sum_{i=1}^n S_iS_i^*$ is
isomorphic to the compact operators $\fK$, and the quotient by this
ideal is $\O_n$.
Let the canonical generators of $\E_n$ be denoted by $\fs_1,\dots,\fs_n$.
Then every such $n$-tuple of isometries arises from a $*$-representation
$\sigma$ of $\E_n$ (write $\sigma \in \Rep(\E_n)$) as $S_i =
\sigma(\fs_i)$.  

We shall call a representation $\sigma$ of $\fA_n$
\textit{$*$-extendible} if $\sigma$ is the restriction to $\fA_n$ of a
$*$-representation of $\E_n$ or $\O_n$ to the canonical copy of
$\fA_n.$   It is easy to see that $\sigma$ is $*$-extendible if and only
if $\sigma(\fS_i)$ are isometries with orthogonal ranges; or
equivalently, $\sigma$ is contractive and $\sigma(\fs_i)$ are
isometries.

Let $\Fn$ denote the unital free semigroup on $n$ letters. (Probably
we should use the algebraist's term `monoid' here, but our habit of
using the term semigroup is well entrenched.)  This semigroup consists
of all words $w$ in $1,2,\dots,n$ including the empty word $\mt$.  The
Fock space $\Hn$ has an orthonormal basis $\{ \xi_w : w \in \Fn \}$,
and is the natural Hilbert space for the left regular representation
$\lambda$.  This representation has generators, denoted by
$L_i:=\lambda(\fs_i)$, which act by $L_i \xi_w = \xi_{iw}$.  The
\wot-closed algebra that they generate is denoted by $\fL_n$.

In general, each $n$-tuple $S_1,\dots,S_n$ will generate a unital
algebra, and the \wot-closure will be denoted by $\fS$.
When a representation $\sigma$ of $\E_n$ is given and $S_i =
\sigma(\fs_i)$, we may write $\fS_\sigma$ for clarity.
For each word $w = i_1\dots i_k$ in $\Fn$, we will use the
notation $S_w$ to denote the corresponding operator $S_{i_1}\cdots
S_{i_k}$.  In particular, $L_w \xi_v = \xi_{wv}$ for $w,v \in \Fn$.

It is a basic fact of C*-algebra theory that every representation of
$\E_n$ splits as a direct sum of the representation induced from its
restriction to $\fK$ and a representation that factors through the
quotient by $\fK$.  However, $\fK$ has a unique irreducible
representation, and it induces the left regular representation
$\lambda$ of $\Fn$, described above.  So $\sigma \simeq
\lambda^{(\alpha)} \oplus \tau$ where $\alpha$ is some cardinal and
$\tau$ is a representation of $\O_n$.   This is equivalent to the
spatial result known as the Wold decomposition.  The Wold
decomposition is the observation that the range $\M$ of the projection
$I - \sum_{i=1}^n S_iS_i^*$ is a wandering subspace, meaning that the
subspaces $\{ S_w \M : w \in \Fn \}$ are pairwise orthogonal, and
together span the subspace $\fS[\M]$.  Any
orthonormal basis for $\M$ will consist of wandering vectors which
generate orthogonal copies of the left regular representation;
moreover, the
restriction of the $S_i$ to $\fS[\M]^\perp$ will be a representation
which factors through $\O_n$.  We call the representation $\tau$ the
\textit{Cuntz part} of $\sigma$, and when $\alpha=0$,
i.e. when $\sum_{i=1}^n \sigma(\fs_i\fs_i^*) =I$, we say simply that
the representation $\sigma$ is of \textit{Cuntz type}.

Recall \cite{DP1} that every $A\in \fL_n$ has
a Fourier series $A \sim \sum_{w \in \bF_n^+} a_w L_w$
determined by $A \xi_\mt = \sum_{w \in \bF_n^+} a_w \xi_w$.
The representation $\lambda$ is a canonical completely isometric map 
from $\fA_n$ into $\fL_n$ which sends $\fs_i$ to $L_i$.    
Hence elements of $\fA_n$ inherit corresponding Fourier series,  and
we will write $A \sim \sum_{w \in \bF_n^+} a_w \fs_w$.
The functional $\phi_0$ reads off the coefficient $a_\mt$.
The kernel of $\phi_0$ in $\fA_n$ and $\fL_n$ are denoted by
$\fA_{n,0}$ and $\fL_{n,0}$ respectively.
These are the norm and \wot-closed ideals, respectively, generated by
the generators $\fs_1,\dots,\fs_n$ and $L_1,\dots,L_n$.
Even when $\phi_0$ is not defined on a free semigroup algebra $\fS$, we
still denote by $\fS_0$ the \wot-closed ideal generated by
$S_1,\dots,S_n$.
This will either be codimension one or equal to the entire algebra.
 
The ideals $\fA_{n,0}^k$ and $\fL_{n,0}^k$ consist of those elements with
zero Fourier coefficients for all words $w$ with $|w| < k$; and are
generated as a right ideal by $\{ \fs_w : |w| = k \}$. 
Moreover \cite{DP2}, each element in $\fL_{n,0}^k$ may be uniquely
represented as $A = \sum_{|w|=k} L_w A_w$ and $\|A\|$ is equal to the
norm of the column operator with entries $A_w$.

One can recover an element of $\fA_n$ or $\fL_n$ from its Fourier
series in the classical way using a summability kernel.
For $t\in\bT$, let $\alpha_t$ be the gauge automorphism of $\O_n$
determined by the mapping $\fs_i\mapsto t \fs_i$.  
Let $V_n(t)=\sum_{k=-2n-1}^{2n+1} c_k t^k$ be the de la Vall\'{e}e
Poussin summability kernel on $\bT$ from harmonic analysis. 
Recall that $V_n$ is a trigonometric polynomial of degree $2n+1$ with
Fourier transform $\hat{V}_n(k)=1$ for $|k|\leq n+1.$
Let $m$ be normalized Lebesgue measure on $\bT$. 
Define linear maps $\Sigma_k$ on $\O_n$ by 
\[
  \Sigma_k(X)= \int_{\bT} V_k(t)\alpha^{-1}_t(X)\, dm(t). 
\]
Then $\Sigma_k$ is a unital completely positive map on $\O_n$ which
leaves $\fA_n$ invariant and moreover, for every $X\in \O_n$,
$\Sigma_k(X)$ converges in norm to $X$.  It has the additional
property that the Fourier coefficients of $\Sigma_k(X)$ agree with
those of $X$ up to the $k$-th level.  
Indeed, if $A \sim \sum_{w \in \bF_n^+} a_w L_w$ lies in $\fA_n$, then
$\Sigma_k(A) = \sum_{|w| \le 2k+1} c_{|w|} a_w \fs_w$.  
Notice that for $A \in \fL_n$, $\Sigma_k(A)$ converges to $A$ in the 
strong operator topology.

Let $\sigma$ be a $*$-extendible representation and let
$\fS=\fS_\sigma$.  We now recall some facts from \cite{DKP} regarding
the ideals $\fS_0^k$.  The intersection $\fJ$ of these ideals is a
left ideal of the von Neumann algebra $\fW$ generated by the $S_i$;
therefore, $\fJ$ has the form $\fW P_\sigma$ for some projection
$P_\sigma \in \fS$. (When the context is clear, we will write $P$
instead of $P_\sigma$.)  The Structure Theorem for free semigroup
algebras \cite{DKP} shows that $P$ is characterized as the largest
projection in $\fS$ such that $P\fS P$ is self-adjoint.  Moreover
$P^\perp \H$ is invariant for $\fS$ and when $P\neq I$, the
restriction of $\fS$ to the range of $P^\perp$ is canonically isomorphic to
$\fL_n$.  Indeed, the map taking $S_i|_{P^\perp\H}$ to $L_i$ extends
to a completely isometric isomorphism which is also a weak-$*$
homeomorphism.  Algebras which are isomorphic to $\fL_n$ are called
\textit{type L}.  When $P\neq I$, the restriction of $\sigma$ to the
range of $P^\perp$ again determines a $*$-extendible representation of
$\fA_n$, and we call this restriction the \textit{type L part} of
$\sigma$.

%%%%%%%%%%%%%%%%%%%%%%%%
\section{Absolute Continuity}\label{S:abs}

In the study of the disk algebra, those functionals which are
absolutely continuous to Lebesgue measure play a special role.
Of course, the Shilov boundary of the disk algebra is the unit
circle, and the Lebesgue probability measure $m$ is Haar measure on it. 
Moreover, every representing measure for evaluation at points interior to
the disk is absolutely continuous.  
We have been seeking an appropriate analogue of this for free semigroup 
algebras for some time. 
That is, which functionals on the non-commutative disk algebra $\fA_n$
should be deemed to be absolutely continuous? 
Unfortunately, there is no clear notion of boundary or representing
measure.  
However there is a natural analogy, and we propose it here.

Our starting point is the left regular representation of the semigroup
$(\bN,+)$.  Under this representation, the generator of $(\bN,+)$ is
mapped to the unilateral shift $S$ and elements of $\bA(\bD)$ are
analytic functions $h(S)$ of the shift, which may be regarded as
multipliers of $H^2(\bT)$.  With this perspective, every vector
functional $h\mapsto \ip{h(S)f_1,f_2}=\int_\bT hf_1\overline{f_2}\,
dm$ corresponds to a measure which is absolutely continuous with
respect to Lebesgue measure.  

On the other hand, suppose $\phi$ is a
functional on $\bA(\bD)$ given by integration over $\bT$ by an
absolutely continuous measure, so that $\phi(h)=\int_\bT hf\, dm$ for
some $f\in L^1(\bT)$.  It is not difficult to show that such
functionals on $\bA(\bD)$ can be approximated by vector functionals
from the Hilbert space of the left regular representation.  Moreover,
if one allows infinite multiplicity, one can represent $\phi$ as a
vector state, that is,  there are vectors
$x_1$ and $x_2$ in $H^{2\,(\infty)}$ such that 
\[
 \phi(h) = \ip{ h(S^{(\infty)}) x_1, x_2 } \qand
 \|\phi\| = \|x_1\| \, \|x_2\| .
\]

Another view is that the absolutely continuous functionals on $\AD$
are the functionals in the predual of $H^\infty(\bT)$.
Our analogue of this algebra is $\fL_n$.  

So we are motivated to make the following definition:

\begin{defn}\label{abs}
For $n\geq 2$, a functional on the non-commutative disk algebra $\fA_n$ is
\textit{absolutely continuous} if it is given by a vector state on $\fL_n$;
i.e. if there are vectors $\zeta, \, \eta \in \ell^2(\bF_n^+)$ so that
$\phi(A) = \ip{ \lambda(A) \zeta, \eta }$.
Let $\fA_n^a$ denote the set of all absolutely continuous
functionals on $\fA_n$.
\end{defn}

For $n\geq 2$, $\fL_n$ has
enough ``infinite multiplicity'' that it is unnecessary to take the
closure of vector functionals; in fact we shall see shortly that
$\fA_n^a$ is already norm closed.

The following result shows that the notion of being representable
as a vector state and being in the predual of $\fL_n$ are
equivalent to each other and to a natural norm condition
on the functional.

\begin{prop}\label{P:abs}
For $\phi \in \fA_n^*$, the following are
equivalent:
\begin{enumerate}
\item $\phi$ is absolutely continuous.
\item $\phi\circ\lambda^{-1}$ extends to a weak-$*$
continuous functional on $\fL_n$.
\item $\dlim_{k\to\infty} \|\phi|_{\fA_{n,0}^k}\| = 0$.
\end{enumerate}
Moreover, given $\eps >0$, the vectors $\zeta$ and $\eta$ may be chosen so that
$\|\zeta\| \, \|\eta\| < \|\phi\| + \ep$.
\end{prop}

\Prf
(1) implies (2) by definition.
The converse follows from \cite{DP1} where it is shown that
every weak-$*$ continuous functional on $\fL_n$ is given by a vector
state.  
The norm condition on the vectors $\zeta$ and $\eta$ is also obtained
there.

Next, suppose (1) holds.  
The map $\lambda$ carries $\fA_{n,0}^k$ into $\fL_{n,0}^k$.
Let $Q_k$ denote the projection of $\Hn$ onto 
$\overline{\spn}\{\xi_w: w\in\bF_n^+, |w|\geq k\}$. 
Elements of $\fL_{n,0}^k$ are characterized by $A = Q_k A$.  
Thus for $A \in \fA_{n,0}^k$,
\[
 |\phi(A)| = | \ip{ \lambda(A) \zeta, \eta} |
           = | \ip{ \lambda(A) \zeta, Q_k \eta} |
 \le \|A\| \, \|\zeta\| \, \| Q_k \eta \| .
\]
Hence it follows that
\[
 \dlim_{k\to\infty} \|\phi|_{\fA_{n,0}^k}\| \le 
 \dlim_{k\to\infty} \|\zeta \| \, \| Q_k \eta \| = 0 .
\]

Conversely suppose that (3) holds.  
Then given $A\in\fA_n$, we use the fact that $\Sigma_k(A)$ converges to
$A$ in norm. 
Note that when $m\geq k$, $\Sigma_k(X)-\Sigma_m(X)$ belongs to
$\fA_{n,0}^k$ and has norm at most $2\|A\|$.
It follows therefore that the adjoint maps satisfy, 
\[
 \| \Sigma_k^*(\phi) - \Sigma_m^*(\phi) \|  \le
 2 \|\phi|_{\fA_{n,0}^k}\| .
\]
So $\Sigma_k^*(\phi)$ is a Cauchy sequence in $\fA^*$. 

We claim that $\Sigma_k^*(\phi)$ is absolutely continuous.
Indeed, consider the Fourier series $A\sim \sum_{w\in \bF_n^+}
a_w\fs_w$.  Then 
\begin{align*}
 \Sigma_k^*(\phi)(A) &= \sum_{|w| \le 2k+1} c_{|w|} a_w \phi(\fs_w)
 \\&=
 \sum_{|w| \le 2k+1} c_{|w|} \phi(\fs_w) \ip{\lambda(A) \xi_\mt,\xi_w}
 \\&=
 \ip{ \lambda(A) \xi_\mt, 
 \tsum_{|w|\le 2k+1}\!\!\! c_{|w|} \phi(\fs_w) \xi_w}.
\end{align*}
{From} (1) implies (2), we found that the set of absolutely continuous functionals is norm
closed.  Hence the limit $\phi$ is also absolutely continuous.
\qed

The following is immediate.

\begin{cor}\label{C:predual}
The set $\fA_n^a$ is the closed subspace of the dual of $\fA_n$
which forms the predual of $\fL_n$.
\end{cor}

\begin{defn}\label{D:abs_cnts_vector}
Let $\sigma$ be a $*$-extendible representation of $\fA_n$ on the
Hilbert space $\H_\sigma$.
A vector $x\in \H_\sigma$ is called an \textit{absolutely continuous}
vector if the corresponding vector state taking $A\in\fA_n$ to
$\ip{\sigma(A)x,x}$ is absolutely continuous. 
\end{defn}

Another straightforward but useful consequence is:

\begin{cor}\label{C:typeL_is_abs}
If $\sigma$ is a $*$-extendible representation of $\fA_n$ and $x,y$ are
vectors lying in the type L part of $\fS = \wotclos{\sigma(\fA_n)}$ 
$($or even in the type L part of 
$\fT = \wotclos{(\sigma\oplus\lambda)(\fA_n)}\ )$, then 
$\phi(A)  = \ip{ \sigma(A) x,y }$ is absolutely continuous. 
In particular, every vector lying in the type L part of $\H_\sigma$ is
absolutely continuous.
\end{cor}

\Prf Relative to $\H\oplus\Hn$, the structure projection
$P_{\sigma\oplus \lambda}$ for $\fT$ decomposes as $P_1\oplus 0$,
with $P_1\leq P_\sigma$.  By considering vectors of the form $x\oplus
0$, where $x$ is in the range of $P_\sigma$, we may regard the type L
part of $\fS$ as contained in the type L part of $\fT$.  Thus, we may
assume to be working with the representation $\sigma\oplus \lambda$
from the start.  By \cite[Theorem~1.6]{DKP}, the type L part of
$\sigma\oplus \lambda$ is spanned by wandering vectors .  For any
wandering vector $w$, the functional $\phi_w(A) = \ip{
(\sigma\oplus\lambda)(A) w, y \oplus 0 }$ is absolutely continuous
because the cyclic subspace $(\sigma\oplus\lambda)(\fA_n)[w]$ is
unitarily equivalent to $\Hn$ and $y \oplus 0$ may be replaced with
its projection into this subspace.  By the previous corollary, the set
of absolutely continuous functionals is a closed subspace.  Taking
linear combinations and limits shows that $\phi$ is in this closure,
and hence also absolutely continuous.  \qed

Now we wish to develop a connection between absolute continuity
and certain intertwining operators.  

\begin{defn}\label{D:intertwiner}
Let $\sigma$ be a $*$-extendible representation of $\fA_n$ on a
Hilbert space $\H$ with generators $S_i = \sigma(\fs_i)$.
Say that an operator $X \in \B(\Hn,\H)$ \textit{intertwines} $S$
and $L$ if $S_iX = X L_i$ for $1 \le i \le n$.
Let $\X(\sigma)$ denote the set of all such intertwiners.  We denote
the range of the subspace $\X(\sigma)$ by $\vac(\sigma)$, that is,
$\vac(\sigma)=\X(\sigma)\ell^2(\bF_n^+)$.  
\end{defn}

Notice that $\vac(\sigma)$ is an invariant linear manifold for
$\sigma(\fA_n)$.  The next result shows that $\vac(\sigma)$ is also
closed and equals the set
of absolutely continuous vectors for $\sigma$.

\begin{thm}\label{intertwiner}
For a representation $\sigma$ of $\E_n$ on $\H$, let  $Q$ be the
structure projection of   
$\fT = \wotclos{(\sigma\oplus\lambda)(\fA_n)}.$  The
following statements hold. 
\begin{enumerate}
\item[i)] For $x,y\in \vac(\sigma)$, the functional
  $\psi(A)=\ip{\sigma(A)x,y}$ is absolutely continuous on $\fA_n$.
\item[ii)] If $x\in\H$ and $\psi(A)=\ip{\sigma(A)x,x}$ is absolutely
  continuous, then $x\in \vac(\sigma).$
\item[iii)] The subspace $\vac(\sigma)$ is closed and is
  $\sigma(\fA_n)$-invariant.   Moreover, it is the 
subspace of $\H$ corresponding to the projection onto $\H$ of the type
L part of $(\sigma\oplus\lambda)(\fA_n)$, that is, $\vac(\sigma)=\ran
(P_\H Q^\perp)$.
\end{enumerate}
\end{thm}

\Prf Let $x,y\in\vac(\sigma)$, and choose 
 vectors $\zeta, \,\eta$ in
$\ell^2(\bF_n)$ and $X,\, Y\in \X(\sigma)$ with $X\zeta = x$ and $Y \eta = y$.
Then 
\[
  \phi(A) = \ip{ \sigma(A) x,y } = \ip{ \sigma(A) X\zeta, Y\eta }
        = \ip{ X \lambda(A) \zeta, Y\eta }
           = \ip{ \lambda(A) \zeta, X^*Y\eta },
\]
so $\phi$ is absolutely continuous.

Suppose that $\psi(A) = \ip{\sigma(A) x, x }$ is absolutely
continuous, say $\psi(A) = \ip{\lambda(A) \zeta, \eta }$.  Theorem~1.6
of \cite{DKP} shows that $x \oplus \zeta$ is a cyclic vector for an
invariant subspace $\M$ of $(\sigma\oplus\lambda)(\fA_n)$ on which the
restriction is unitarily equivalent to $\lambda$.  Indeed, while the
hypothesis of that theorem requires that $\sigma$ be type L, this
condition is used only to establish that $\psi$ is absolutely
continuous (in our new terminology).  It is evident that a subspace of
this type is the range of an intertwining isometry $V \in
\X(\sigma\oplus\lambda)$.  Let $X = P_\H V$.  Then $X$ intertwines $S$
and $L$.  Moreover, since $x \oplus \zeta$ is in the range of $V$, it
follows that $x$ is in the range of $X$, so (ii) holds.

We now push this argument a little further.  Observe that as in the
 proof of \cite[Theorem~1.6]{DKP}, given $t>0$, $\zeta$ may be
 replaced by $t\zeta$.
 Therefore, if $x\in\vac(\sigma)$, the argument of the previous
 paragraph also shows that $x$ belongs to the closed span of the
 wandering vectors for $(\sigma\oplus\lambda)(\fA_n)$.  Thus $x$
 belongs to the type L part of $(\sigma\oplus\lambda)(\fA_n)$, whence
 $\vac(\sigma)\subseteq \ran(P_\H Q^\perp).$ Conversely, since $P_\H$
 and $Q$ commute, any vector $x\in\ran(P_\H Q^\perp)$ lies in the type
 L part of $\fT$, and thus $\psi(A)=\ip{Ax,x}$ is absolutely
 continuous by Corollary~\ref{C:typeL_is_abs}. But then
 $x\in\vac(\sigma)$ by part (ii).  So $\ran(P_\H
 Q^\perp)=\vac(\sigma)$.  That $\vac(\sigma)$ is closed is now
 obvious.  \qed

We now give a condition sufficient for the existence of wandering
vectors. 

\begin{thm}\label{T:inter_wander}
Let $X$ belong to  $\X(\sigma)$.
Then the following statements are equivalent.
\begin{enumerate}
\item[i)] The representations  $\sigma|_{\ol{\ran X}}$ and $\lambda$
  are  unitarily equivalent;
\item[ii)] $\ol{\ran X} = \fS[w]$ for some wandering vector $w$;
\item[iii)] $X^*X = R^*R$ for some non-zero $R\in \fR_n = \fL_n'$.
\end{enumerate}
In particular, this holds if $X$ is bounded below.
\end{thm}

\Prf The equivalence of (i) and (ii) is clear from the
 definitions.

To obtain (iii) $\Rightarrow$ (i), suppose that $X^*X = R^*R$.
By restricting $\sigma$ to the invariant subspace $\ol{\ran X}$, we
may suppose that $X$ has dense range, and that $X \xi_\mt$ is a
cyclic vector. 
We now show that $\sigma$ is  equivalent to
$\lambda$.
  
Since $R\in\fR_n$, \cite[Corollary~2.2]{DP1}, shows that $R$ factors as the
product of an isometry and an outer operator in $\fR_n$. The equality
$X^*X=R^*R$ is unchanged if the isometry is removed, so we may assume
that $R$ has dense range.  Since $X$ and $R$ have the same positive
part,  there is an
isometry $V$ such that $X = VR$ and $\ran V = \ol{\ran X}$; whence $V$
is unitary.  Then
\[
  (S_iV-VL_i)R = S_i(VR) - (VR)L_i = 0.
\]
Therefore $V$ intertwines $S$ and $L$ and so $\sigma|_{\ran V}$ is
equivalent to $\lambda$.

Finally, we show (ii) $\Rightarrow$ (iii).  If there is an isometry $V
\in \X(\sigma)$ with $\ran V = \ol{\ran X}$, then by again restricting
to this range, we may assume that $V$ is unitary, so that $\sigma$ is
equivalent to $\lambda$.  So $V^*S_i = L_i V^*$.  Hence
\[
 (V^*X)L_i = V^* S_i X = L_i (V^* X) ;
\]
whence $R := V^*X$
belongs to $\fL_n' = \fR_n$.
Therefore $X = VR$ and so $X^*X = R^*R$.

Now suppose that $X$ is bounded below.
Again we may suppose that $X$ has dense range, hence $X$ is
invertible.

Consider the Wold decomposition of $S$.
The Cuntz part is supported on 
\begin{align}
  \N &:= \bigcap_{k\ge1} \sum_{|w|=k} \ran S_w 
       = \bigcap_{k\ge1} \sum_{|w|=k} S_w X \Hn \label{e:N}
  \\&= \bigcap_{k\ge1} \sum_{|w|=k} X L_w \Hn 
     = X \bigcap_{k\ge1} \sum_{|w|=k} L_w \Hn = \{0\} \notag
\end{align}
Hence $\sigma$ is a multiple of $\lambda$. 
Since $\ran X$ has a cyclic vector $X \xi_\mt$, $\sigma$
has multiplicity one, and thus is equivalent to $\lambda$.
\qed  

As an immediate corollary, we note the existence of wandering
vectors is characterized by a structural property of $\X(\sigma)$.

\begin{cor}\label{wandcondition}
Let $\sigma$ be a representation of $\E_n$ on $\H$ with generators
$S_i=\sigma(\fs_i)$.  Then $S$ has a wandering vector if
and only if there exists $X\in\X(\sigma)$ such that $X$ is bounded
below.
\end{cor}

\Prf If $\eta\in\H$ is wandering for $\fS$, then the isometric map
determined by $X\xi_w = w(S)\eta$ belongs to $\X(\sigma)$.  The converse
follows from the theorem.
\qed

\begin{rem}
If one only has $X^*X \ge R^*R$ for a non-zero $R\in \fR_n$, one may
still deduce that $\ran X$ has wandering vectors.  
To do this, use Douglas' Lemma \cite{Do} to factor $R = YX$. Then argue as
in Theorem~\ref{T:inter_wander} that $YS_i = L_iY$.
Then with $\N$ as in \eqref{e:N}, one can show that $Y \N = \{0\}$.
Since $Y$ has dense range, $\sigma$ has a summand equivalent to
$\lambda$.
Moreover, since the range of an intertwiner consists of absolutely
continuous vectors, the existence of this summand and
Lemma~\ref{spanning} below show that the range of $X$ is spanned by
wandering vectors.
\end{rem} 

\begin{eg}
There are intertwining maps whose range is not equivalent to $\lambda$.
For example, consider the atomic representation of type
$\pi_{z_2^\infty}$ \cite[Example~3.2]{DP1}.
Then the restriction of $S_2$ to the \textit{spine} 
$\ell^2(\bZ \times \{0\})$ is the bilateral shift.

Observe that there is a summable sequence $(a_k)_{k\in\bZ}$ such
that $\sum_{k\in\bZ} a_k \xi_{k,0}$ is cyclic for the bilateral shift.
Indeed, Beurling's Theorem states that the (cyclic) invariant
subspaces of the bilateral shift, considered as $M_z$ on $L^2(\bT)$,
have the form $L^2(E)$ for a measurable subset $E$ of $\bT$ or the
form $w H^2$ where $|w| = 1$ a.e.  Thus if a function $g$ vanishes on
a set of positive measure, it generates $L^2(\operatorname{supp}(g))$.
On the other hand, if there is an outer function $f$ in $H^2$ with
$|f|=|g|$ a.e., then the cyclic subspace is $w H^2$ where $w = g/f$.
This occurs if and only if $\log |g|$ belongs to $L^1(\bT)$.  So
choose a $C^2$ function $g$ on $\bT$ which vanishes at a single point
in such a way that $\log |g|$ is not integrable.  For example, make
$g(\theta) = e^{-1/|\theta|}$ near $\theta = 0$ and smooth.  Lying in
$C^2$ guarantees that the Fourier coefficients are summable.

For each $k \in \bZ$, there is an intertwining isometry $V_k$ with 
$V_k \xi_\mt = \xi_{k,0}$.
Then $V = \sum_{k\in\bZ} a_k V_k$ is an intertwiner.
Moreover, $V \xi_\mt$ is cyclic for this Cuntz representation.
So $V$ has dense range; but the representation $\pi_{z_2^\infty}$ is not
equivalent to $\lambda$.
\end{eg}

\begin{rem}\label{poisson}
Consider the completely positive map on $\B(\H)$ given by
$\Phi(A) = \sum_i S_i A S_i^*$. 
Suppose that $X$ intertwines $S$ and $L$. 
Then 
\[
 \Phi^k(XX^*) = \sum_{|w|=k} S_w X X^* S_w^* = 
 X \sum_{|w|=k} L_w L_w^* X^* = X Q_k X^* \le XX^* .
\] 
Moreover, $\sotlim_k \Phi^k(X^*X) = 0$.  
This latter condition is called \textit{purity} by Popescu
\cite{Pop_poisson}. 
Under these two hypotheses, namely $\Phi(D) \le D$ and $\sotlim_k
\Phi^k(D) = 0$, Popescu proves the converse, that  $D = X^*X$ for an
intertwiner $SX = XL^{(\infty)}$ using his Poisson transform.
\end{rem}

%%%%%%%%%%%%%%%%%%%%%%%%%%%%%%%%%%%%%%%%%%%%%%%%%%%%%%%%%%%%
\section{Wandering vectors and absolute continuity}\label{S:wander}

In \cite{DKP}, we showed that in the presence of summands which
contain wandering vectors, the entire type L part is spanned by
wandering vectors.
In this section, we use the ideas of the previous section to
strengthen this significantly by showing that the presence of one
wandering vector implies that the type L part is spanned by wandering
vectors.
We then consider the various ways in which a representation can
appear to be type L.

\begin{defn}\label{typeLclones}
Let $\sigma$ be a $*$-extendible representation of $\fA_n$. 
We say that $\sigma$ is \textit{type L} if the free semigroup algebra
generated by $\sigma(\fs_1)$, \dots, $\sigma(\fs_n)$ is type L.

A representation $\sigma$ is \textit{weak type L} if 
$\sigma \oplus \lambda$ is type L.  

A representation $\sigma$ is \textit{weak-$*$ type L} if 
$\sigma^{(\infty)}$ is type L.  
 
The representation $\sigma$ of $\fA_n$ is \textit{absolutely
continuous} if every vector state $\psi(A) = \ip{ \sigma(A) x,x}$ is
absolutely continuous.

\end{defn}

Notice that the restriction of a $*$-extendible representation
$\sigma$ of $\fA_n$ to the
invariant subspace $\vac(\sigma)$ produces an absolutely continuous
representation.  We call this restriction the \textit{absolutely
  continuous part} of $\sigma$.

\begin{lem}\label{spanning}
If $\sigma$ is absolutely continuous and has a wandering vector, 
then $\H$ is spanned by its wandering vectors.
In particular, $\sigma$ is type L.
\end{lem}

\Prf 
Let $\eta$ be a wandering vector in $\H$, and set $\H_0 = \fS[\eta]$.
Let $V$ be the isometry in $\X(\sigma)$ mapping $\Hn$ onto $\H_0$.

By Theorem~\ref{intertwiner}, every vector $x\in\H$ is in the range
of some intertwining map $X \in \X(\sigma)$.
We may assume that $\|X\|=1/2$.
Then $V \pm X$ are intertwiners which are bounded below.
By Theorem~\ref{T:inter_wander}, the ranges of these two intertwiners
are the ranges of isometric intertwiners, and thus are spanned by
wandering vectors.
But the range of $X$ is contained in the sum of the ranges of
$V \pm X$; and hence $x$ is contained in the span of all wandering
vectors.
\qed 

\begin{cor}\label{span_wander}
If $\sigma$ is any representation of $\E_n$ such that $\sigma(\fA_n)$
has a wandering vector, then the span of the wandering vectors for
$\sigma(\fA_n)$ is $\vac(\sigma)$.
\end{cor}

\Prf Any wandering vector is an absolutely continuous vector, so 
simply restrict $\sigma$ to the $\sigma(\fA_n)$-invariant subspace
consisting of absolutely continuous vectors and apply the lemma.
\qed

We now delineate the various type L forms, and their relationships as we
know today.  There are no known examples of absolutely continuous
representations without wandering vectors.

\begin{thm}\label{T:equiv_ac}
Consider the following conditions for a 
$*$-extendible representation $\sigma$ of $\fA_n$:
\begin{enumerate}
\item[\textrm{(1a)}] $\sigma$ is absolutely continuous
\item[\textrm{(1b)}] $\sigma \oplus \lambda$ is type L $($i.e. $\sigma$ is weak
                type L$)$
\item[\textrm{(1c)}] $\sigma \oplus \tau$ is type L for any $($all\/$)$
type L representation $\tau$.
\item[\textrm{(2a)}] $\sigma^{(\infty)}$ is type L $($i.e. $\sigma$ is
weak-$*$ type L$)$
\item[\textrm{(2b)}] $\sigma$ is absolutely continuous and 
$\sigma^{(\infty)}$ has a wandering vector
\item[\textrm{(2c)}] $\sigma^{(\infty)}$ is spanned by wandering vectors
\item[\textrm{(3a)}] $\sigma$ is type L 
\item[\textrm{(3b)}] $\sigma$ is absolutely continuous and 
$\sigma^{(k)}$ has a wandering vector for some
finite $k$
\item[\textrm{(3c)}] $\sigma^{(k)}$ is spanned by wandering vectors  for
some finite $k$
\item[\textrm{(4a)}] $\sigma$ is absolutely continuous and 
has a wandering vector
\item[\textrm{(4b)}] $\sigma$ is spanned by wandering vectors
\end{enumerate}
Then properties with the same numeral are equivalent, 
and larger numbers imply smaller. 
\end{thm}

\Prf
(1a) $\Rightarrow$ (1b).  
If $\sigma$ is absolutely continuous, then $\sigma\oplus\lambda$ is
absolutely continuous and has a wandering vector.
Thus by Lemma~\ref{spanning}, $\sigma\oplus \lambda$ is spanned by its
wandering vectors, and so is type L.

(1b) $\Rightarrow$ (1a). Since $\sigma\oplus \lambda$ is type L and has a
wandering vector, Lemma~\ref{spanning} shows that
$\sigma\oplus\lambda$ is spanned by its wandering vectors.  
Thus $\sigma\oplus\lambda$ is absolutely continuous, and hence so is
$\sigma.$   

(1a) $\Rightarrow$ (1c): If $\tau$ is any type L representation, there
is an integer $p$ so that $\tau^{(p)}$ has a wandering vector.  Thus
$(\sigma\oplus\tau)^{(p)}$ is absolutely continuous and has a
wandering vector, and so is also type L.  However being type L is not
affected by finite ampliations, as this has no effect on the
\wot-closure.  So $\sigma \oplus \tau$ is type L.  

(1c) $\Rightarrow$ (1a): If $\sigma
\oplus \tau$ is type L for some type L representation $\tau$, then
$\sigma\oplus\tau$  is absolutely continuous. 
By considering  vectors of the form $x\oplus 0$, we find that
$\sigma$ is absolutely continuous.  
So (1a), (1b), and (1c) are all equivalent.

If $\sigma$ is weak-$*$ type L, then $\sigma^{(\infty)}$ has a finite
ampliation which is spanned by wandering vectors.  But of course this
ampliation is equivalent to $\sigma^{(\infty)}$, so (2a) implies
(2c). Clearly, if (2c) holds, then every vector in $\H^{(\infty)}$ is
  absolutely continuous, so in particular, $\sigma$ is an absolutely
  continuous representation; thus (2b) holds.  If $\sigma$ is absolutely
  continuous and
 $\sigma^{(\infty)}$ has a wandering vector, then $\H^{(\infty)}$ is
spanned by wandering vectors and thus $\sigma^{(\infty)}$ is type L.    So
(2a), (2b), and (2c) are all equivalent and imply (1).

The equivalence of (3a), (3b) and (3c) follows from
\cite[Corollary~1.9]{DKP}, and evidently implies (2).

By Lemma~\ref{spanning}, (4a) and (4b) are equivalent and clearly imply
(3).
\qed

It is worthwhile examining the various weaker notions of type L in
light of the Structure Theorem for Free Semigroup Algebras \cite{DKP}.
Let $\sigma$ be a representation of $\E_n$ and let $\fS$ and $\fW$ denote the
corresponding free semigroup algebra and von Neumann algebra respectively.
Then there is a projection $P$ in $\fS$ characterized as the largest projection
in $\fS$ for which $P \fS P$ is self-adjoint.
Then $\fS = \fW P + \fS P^\perp$, $P^\perp \H$ is invariant for $\fS$
and $\fS P^\perp$ is type L.  
We wish to break this down a bit more.

\begin{defn}\label{type}
A representation $\sigma$ of $\E_n$ or $\O_n$ is of \textit{von
Neumann type} if the corresponding free semigroup algebra $\fS$ is a
von Neumann algebra.  If $\sigma$ has no summand of either type L or
von Neumann type, say that it is of \textit{dilation type}.  We also
will say that $\sigma$ is weak-$*$ of some type if $\sigma^{(\infty)}$
is of that type.
\end{defn}

A very recent result of Charles Read \cite{Read} shows that there can
indeed be representations of von Neumann type.

The reason for the nomenclature \textit{dilation type} is that after
all summands of von Neumann type and type L are removed, the remainder
must have a non-zero projection $P$ prescribed by the structure
theorem such that $P\H$ is cyclic and $P^\perp\H$ is cyclic for
$\fS^*$.  For these algebras, the type L corner must be a multiple of
$\lambda$.  To see this, consider the subspace $\W = \big( \sum_i S_i
P \H \big) \ominus P\H$.  This is a wandering subspace for the type L
part.  It is necessarily non-zero, for otherwise $\fS$ would be a von
Neumann algebra.  Moreover, $\W$ is cyclic for the type L corner because
of the cyclicity of $P\H$. Hence the type L part is equivalent to
$\lambda^{(\Dim \W)}$. This is an observation that was, unfortunately,
overlooked in \cite{DKP}. Hence one sees that the compressions $A_i =
P S_i |_{P\H}$ form a row contraction with $S_i$ as their minimal
isometric dilation (in the sense of Frahzo--Bunce--Popescu).  We
record the most useful part of this for future reference.

\begin{prop}\label{P:dilation}
If $\sigma$ is dilation type, then it has wandering vectors.
In particular, dilation type and weak-$*$ dilation type coincide.
\end{prop}

\Prf The first statement was proven in the preamble.  Once one has
a wandering vector, the span of the wandering vectors includes all
of the absolutely continuous vectors, which includes the weak-$*$
type L part.
\qed

We can now clarify the exceptional case in which there may be pathology.

\begin{prop}\label{P:exception}
Let $\sigma$ be a $*$-extendible representation of $\fA_n$.
If the type L and absolutely continuous parts do not coincide, then
$\sigma$ is of von Neumann type, and decomposes as 
$\sigma \simeq \sigma_a \oplus \sigma_s$ where $\sigma_a$ is absolutely
continuous and $\sigma_s$ has no absolutely continuous part.
\end{prop}

\Prf Decompose $\sigma \simeq \sigma_v \oplus \sigma_d \oplus \sigma_l$
into its von Neumann, dilation and type L parts.
By Proposition~\ref{P:dilation}, if there is a dilation part, then there
are wandering vectors.
So by Corollary~\ref{span_wander}, the type L and absolutely continuous
parts coincide.
Likewise if there is a type L part, the equivalence of (1a) and (1b) in
Theorem~\ref{T:equiv_ac} shows that the type L and absolutely continuous
parts will coincide.
So $\sigma$ is necessarily of von Neumann type.

Since $\vac(\sigma)$ is invariant for $\fS_\sigma$, and $\fS_\sigma$
is a von Neumann algebra, $\vac(\sigma)$ is a reducing subspace for
$\fS_\sigma$.  This gives the desired decomposition  $\sigma\simeq
\sigma_a\oplus \sigma_s$.
\qed

\begin{defn}
Call a $*$-extendible representation $\sigma$ of $\fA_n$ \textit{regular}
if the absolutely continuous and type L parts of $\sigma$ coincide. 
\end{defn}  

\begin{rem}\label{R:pathology}
Proposition~\ref{P:exception} shows that the only pathology that can
occur in the various weak type L possibilities is due to a lack of
wandering vectors.

It is conceivable that a representation is type L but has no
wandering vectors.  Such an algebra is reductive and nonselfadjoint.  
There is no operator algebra known to have this property.
So the (unlikely) existence of such an algebra would yield a counterexample
to a well-known variant of the invariant subspace problem.

 A $*$-extendible representation $\sigma$ which is weak-$*$ type L but
not type L must be von Neumann type by the preceding proposition.
But then $\wsclos{\sigma(\fA_n)}$ would be a weak-$*$ closed
subalgebra isomorphic to $\fL_n$ which is \wot-dense in a von Neumann
algebra.  We have no free semigroup algebra example of this type of
behaviour.  However, Loebl and Muhly \cite{LM} have constructed an operator
algebra which is weak-$*$ closed and nonselfadjoint, but with the
\wot-closure equal to a von Neumann algebra. Therefore it is conceivable that
such a free semigroup algebra could exist.

Finally, one could imagine that $\sigma$ is of weak-$*$ von Neumann type
but absolutely continuous.

Clearing up the question of whether any of these possibilities can
actually occur remains one of the central questions in the area.  
We conjecture that every representation is regular.
Indeed, we would go further and speculate that type L representations
always have wandering vectors.
\end{rem}

%%%%%%%%%%%%%%%%%%%%%%%%%%%%%%%%%%%%%%%%%%%%%%%%%%%%%%%%%%%%
\section{Reflexivity and hyper-reflexivity}\label{S:refl}

In this  section, we establish two reflexivity results that 
extend previous work in light of the previous section.

\begin{thm}
If $\fS$ is a free semigroup algebra which has a wandering vector, then it
is reflexive.
\end{thm}

\Prf By \cite[Proposition~5.3]{DKP}, $\fS$ is reflexive if and only if
the restriction to its type L part is reflexive. Thus, without loss of
generality, we may assume that $\fS$ is type L.  Since $\fS$ is type L
and has a wandering vector, Lemma~\ref{spanning} shows that $\H$ is
spanned by wandering vectors.  Let $W\subseteq\H$ be the set of
all unit wandering vectors. For each $\alpha\in W$, let
$\H_\alpha=\fS[\alpha]$ and let $V_\alpha:\Hn\rightarrow \H_L$ be the
intertwining isometry which sends $\xi_w$ to $S_w\alpha$. Then the
invariant subspaces $\H_\alpha$ span $\H$ and each restriction
$\fS|_{\H_\alpha}$ is unitarily equivalent to $\fL_n$ via $V_\alpha$.

If $T \in \Alg \Lat \fS$, then $\H_\alpha$ is invariant for $T$.
Since $\fL_n$ is reflexive, there is an element $B_\alpha \in \fL_n$
so that $T|_{\H_\alpha} = V_\alpha B_\alpha V_\alpha^*$. For each
$\alpha\in W$, there is an element $A_\alpha \in \fS$ so that
$A_\alpha|_{\H_\alpha} = V_\alpha B_\alpha V_\alpha^*$.  Fix an
element $\alpha_0\in W$, let $V_0=V_{\alpha_0}$ and
$A_0=A_{\alpha_0}.$ We shall show that  $T=A_0$.  By
replacing $T$ with $T-A_0$, we may assume that
$T|_{\H_0}=0$, so that our task is to show $T=0$.  

Given $\alpha\in W$, the operator $X = V_0 + .5 V_\alpha$ is an intertwining
map between $S$ and $L$ which is bounded below.  Moreover, $\M:=\ran X$ is
closed and invariant for $\fS$; hence $\M$ is also invariant for $T$.  But
\[
 T X \xi_\mt = T V_0 \xi_\mt + .5 T V_\alpha \xi_\mt 
 = .5A_\alpha V_\alpha \xi_\mt =: y
\]
belongs to $\H_\alpha \cap \M$.  This implies that there is a vector
$\zeta\in \Hn$ so that $y = X\zeta = V_0 \zeta + .5 V_\alpha \zeta$
belongs to $\H_\alpha$, and thus $V_0 \zeta$ lies in $\H_0 \cap
\H_\alpha$.  If $\zeta=0$ then $y=0$, so that $A_\alpha$ has the
non-zero vector $V_\alpha \xi_\mt$ in its kernel.  Otherwise $V_0
\zeta$ is a non-zero vector in $\H_0 \cap \H_\alpha$ and $A_\alpha
V_0\zeta=TV_0\zeta=0$.  Therefore, $A_\alpha|_{\H_\alpha}$ has
non-trivial kernel.  Hence $B_\alpha$ is an element of $\fL_n$ with
non-trivial kernel.  Since non-zero elements of $\fL_n$ are injective
\cite[Theorem~1.7]{DP1}, we deduce that $B_\alpha=0$.  Hence
$0=A_\alpha|_{\H_\alpha}=T|_{\H_\alpha}$.  Since $\bigvee_{\alpha\in
  W} H_\alpha=\H$, we conclude that $T=0$ as desired.  \qed

Recall that an operator algebra $\fA$ is hyper-reflexive if there is a
constant $C$ so that 
\[
 \dist(T,\fA) \le C \beta_\fA(T) 
 := C \sup_{P \in \Lat \fA} \| P^\perp T P \| .
\]
The known families of hyper-reflexive algebras are fairly small.  It
includes nest algebras \cite{Arv_nest} with constant 1, the analytic
Toeplitz algebra \cite{D_toeplitz} and the free semigroup algebras
$\fL_n$ \cite{DP1}.  Bercovici \cite{Ber} obtained distance constant 3
for all algebras having property $\X_{0,1}$ and also showed that an
operator algebra $\fA$ has property $\X_{0,1}$ whenever its commutant
contains two isometries with orthogonal ranges.  In particular,
$\fL_n$ has property $\X_{0,1}$ when $n\geq 2.$  Bercovici's results
significantly increased the known class of hyper-reflexive algebras.

There is a long-standing open question about whether all von Neumann
algebras are hyper-reflexive, which is equivalent to whether every 
derivation is inner \cite{Chr}.
The missing cases are von Neumann algebras whose commutant are
certain intractable type $II_1$ algebras.
This could include certain type $II_\infty$ representations of $\O_n$,
and hence would apply in our context.
So for the next result, we restrict ourselves to the type L case.

\begin{thm}\label{hyp}
If $\fS$ is a type L free semigroup algebra which has a
wandering vector, then $\fS$ is hyper-reflexive.
\end{thm}

Before giving the proof, we pause for the following remark.

\begin{rem}\label{subhyper}
If $\fS$ is type L and has a wandering vector, then by \cite{DP1} it
has property $\bA_1$ and by \cite{DP3} it even has property
$\bA_{\aleph_0}$.  In particular, Theorem~\ref{hyp} together with a
result from \cite{HN} implies that every
weak-$*$ closed subspace of a type L free semigroup algebra with a
wandering vector is also hyper-reflexive.  Even though
$\X_{0,1}$ is only a bit stronger than $\bA_{\aleph_0}$, we were unable
to show that $\fS$ has it.  So we are unable to apply Bercovici's argument.
Thus, the proof which follows uses methods reminiscent of those used in
\cite{DP1}.

If $\fS$ is type L and has no wandering vector, then as noted in
Remark~\ref{R:pathology}, the algebra will be nonselfadjoint and
reductive.  In particular, it is not reflexive.
\end{rem}

\Prf  Let $T \in \B(\H)$, and set 
 $\beta(T) = \sup_{P\in \Lat \fS} \| P^\perp T P \|$.
Let $x_0$ be a wandering vector of $\fS$.
Then $\fS|_{\fS[x_0]} \simeq \fL_n$. 
Since $\fL_n$ is hyper-reflexive with constant 3, there exists an $A\in
\fS$ with $\| (T- A)|_{\fS[x_0]}|| \le 3 \beta(T)$.
By replacing $T$ with $T-A$, we can assume that
$\| T |_{\fS[x_0]} \| \le 3 \beta(T)$.

Let $x$ be a wandering vector with $x \ne x_0$ and let $V$ be the
isometric intertwiner from $\ell^2(\bF^+_n)$ onto $\fS[x]$ satisfying
$V \xi_w = S_w x$. We shall show that \begin{equation}\label{E:Xbound}
\norm{T|_{\fS[x]}}\leq 26\beta(T).
\end{equation}

Let $x_i = S_i x_0$, for $i =1, 2$.
For $i=0,1,2$, define isometric intertwiners $V_i$ from
$\ell^2(\bF_n^+)$ to
$\H$ by $V_i \xi_w = S_w x_i$ for $w \in \bF^+_n$.

For $i=1,2$, set $T_i = V_i + r V$ where $0 < r < 1/\sqrt2$, 
and define $\N_i = \ran T_i$.
We claim that $\N_1$ and $\N_2$ are at a positive angle to each other;
so that $\N_1 \cap \N_2 =\{0\}$ and $\N_1 + \N_2$ is closed.
Indeed, using $\delta := 1 - r \sqrt2 > 0$,
\begin{align*}
  \| T_1 \xi - T_2 \eta \| &\ge
  \| V_1\xi - V_2\eta \| - r \| V (\xi - \eta) \| \\ &\ge
  \| \xi \oplus \eta \| - r ( \| \xi \| + \| \eta \| ) \;\;\ge
  \delta \| \xi \oplus \eta \| .
\end{align*}
So the natural map of $\N_1\oplus \N_2$ onto $\N_1+\N_2$ is an isomorphism.
 
Observe next that for any
$w\in\bF_n^+$, we have 
\begin{align*}
  \sum_{|w|=k} S_w (\N_1 + \N_2)
  &= \sum_{|w|=k} S_w T_1 \Hn + S_w T_2 \Hn \\
  &= T_1 \sum_{|w|=k} L_w \Hn + T_2 \sum_{|w|=k} L_w \Hn.
\end{align*}
Therefore,
\[
 \lim_{k\to\infty} \sum_{|w|=k} S_w (\N_1+\N_2) =0 .
\]
  As $\sum_{j=1}^{n} S_j \N_i$
has co-dimension one in $\N_i$, we find that.  
\[
 \Dim \Big( \N_1 +\N_2  - \sum_{j=1}^{n} S_j( \N_1 + \N_2) \Big) = 2.
\]
By the Wold decomposition, we deduce that 
 $\fS|_{\N_1+\N_2} \simeq \fL_n^{(2)}$.
This algebra is hyper-reflexive with distance constant 3.
So there is an element $A \in \fS$ such that
$\| (T - A)|_{ \N_1+\N_2 } \| \le 3 \beta(T)$.

Note that $\M := \fS[x_1-x_2] = \ran (T_1 - T_2) \subset \N_1 +
\N_2$; also, $\M$  is a cyclic subspace of $\fS[x_0]$. 
Since $\| T|_\M \| \le 3 \beta(T)$, $\| A|_\M \| \le 6\beta(T)$.  
As $\fS$ is type L, $\| A \| = \| A|_\M \|$, so
$\| T|_{\N_1 + \N_2} \| \le 9 \beta(T).$
We now improve this to an estimate of $\| T|_{\fS[x]} \|$.

Suppose that $y$ is a unit vector in $\fS[x]$. 
Observe that 
\[
  T_1(V^* y) = V_1V^* y + ry
\]
lies in  $\N_1 \subset \N_1 +\N_2$. 
So $\|T (V_1V^* y + ry)|| \leq 9(1+r)  \beta(T)$. 
As $V_1V^* y$ is a unit vector in $\fS[x_0]$,
$\| T V_1V^* y \| \le 3 \beta(T)$. 
Hence 
\[
 \| T y \| \le r^{-1}(12+9r) \beta(T) .
\]
Choosing $r$ sufficiently close to $1/\sqrt2$ 
yields that $\|T|_{\fS[x]} \| \leq 26 \beta(T)$, so \eqref{E:Xbound} holds.

We now can estimate $\norm{T}$.  Fix any unit vector $y \in \H$, and
let $\fT$ be the free semigroup algebra generated by $S_i\oplus L_i$.
Since $\fS$ is type L, by \cite[Theorem 1.6]{DKP} there is a vector
$\zeta \in \ell^2(\bF_n^+)$ with $\|\zeta \|\ < \ep$ such that $\fT[y
\oplus \zeta]$ is a subspace of $\H \oplus \ell^2(\bF_n^+)$ which is
generated by a wandering vector.  Hence $\fT[y \oplus \zeta]$ is the
range of an isometry $W'$ from $\ell^2(\bF_n^+)$ to $\H \oplus
\ell^2(\bF_n^+)$ intertwining $L_i$ with $S_i \oplus L_i$.  Then
$W'':=P_\H W'$ is a contraction in $\B(\Hn,\H)$ satisfying
$S_iW''=W''L_i$.  Moreover, there is a vector $\xi\in\Hn$ of norm $(1
+ \ep^2)^{1/2}$ such that $W''\xi =y$.  Identify $\fS[x_0]$ with
$\ell^2(\bF_n^+)$ via the isometry $V_0\in\B(\Hn,\H)$, and set $W:=
W''V_0^*\in\B(\fS[x_0],\H)$ and $w:= V_0\xi$.

Let $J$ be the inclusion map of $\fS[x_0]$ into $\H$.
For $|t| <1$, consider $V_t = J + t W$.
This is an intertwining map which is bounded below, and thus by
Theorem~\ref{T:inter_wander}, there is a wandering vector $x_t$ of
$\fS$ so that $\ran(V_t) = \fS[x_t]$.  
So 
\[
 \| T (w + t y) \| \le 26 \beta(T) \| w + ty \| . 
\]
Since $\| Tw \| \le 3\beta(T) \|w\|$, if we let $t$ increase to $1$ and
$\ep$ decrease to $0$, we obtain $\| Ty \| \le 55 \beta(T)$.  So
$\norm{T}\leq 55\beta(T)$.  Thus, $\fS$ is hyper-reflexive with constant
at most 55.
\qed

The following proposition is complementary to \cite[Proposition~2.10]{DKP}
showing that if $\fS$ is of Cuntz type, then $\fS'' = \fW$ is a von
Neumann algebra.

\begin{prop}\label{commutant}
Let $\fS$ be a free semigroup algebra acting on a Hilbert
space $\H$ which is not of Cuntz type.
Then $\fS ^{\prime \prime} = \fS$.
\end{prop}

\Prf
Since $\fS$ is not Cuntz type, by the Wold decomposition, it
has a direct summand equivalent to $\fL_n$.
That is, we may decompose the generators $S_1,...,S_n$ as
$S_i = T_i \oplus L_i$ on $\H = H_1 \oplus \ell^2(\bF^+_n)$. 

Let $\fW$ be the von Neumann algebra generated by $\fS$.
By the Structure Theorem \cite[Theorem~2.6]{DKP}, there is a largest
projection $P$ in $\fS$ such that $P \mathfrak S P$ is self-adjoint and
$\fS = \fW P + P^\perp \fS P^\perp$.
Now $\fS'' \subset \fW''= \fW$, so 
$\fS P \subset \fS'' P \subset \fW P = \fS P$;
whence $\fS'' P = \fS P$.

By Theorem~\ref{spanning}, $P^\perp \H$ is spanned by 
wandering vectors.
For any wandering vector $x_\alpha$, let $V_\alpha$ be the canonical
intertwining isometry from $\ell^2(\bF^+_n)$ into $\H$ defined by 
$V_\alpha \xi_w = S_w x_\alpha$ for $w\in\bF^+_n$.
If we select $x_0 = 0 \oplus \xi_\mt$, then $V_0$ maps onto the
free summand.  
It is easy to check that $V_\alpha V_0^*$ commutes with $\fS$.

Let $A \in \fS''$.  Then since $0 \oplus I$ commutes with $\fS$,
$A$ must have the form $A = A_1 \oplus A_2$.
Moreover, $A_2 \in \fL_n'' = \fL_n$ by \cite{DP1}.
There is an element $B \in \fS$ such that $B = B_1 \oplus A_2$.
Subtracting this from $A$, we may suppose that $A = A_1 \oplus 0$.
Then 
\[
 A x_\alpha = A (V_\alpha V_0^*) x_0 = (V_\alpha V_0^*) A x_0 = 0.
\]
Thus $AP^\perp = 0$. 
As above, $AP$ lies in $\fS$, whence $A$ belongs to $\fS$.
\qed

%%%%%%%%%%%%%%%%%%%%%%%%%%%%%%%%%%%%%%%%%%%%%%%%%%%%%%%%%%%%%%%%%%
\section{A Kaplansky Density Theorem}\label{S:kaplansky}

Kaplansky's famous density theorem states that if $\sigma$ is a
$*$-repre\-sent\-ation of a C*-algebra $\fA$, then the unit ball of
$\sigma(\fA)$ is \wot-dense in the ball of the von Neumann algebra
$\fW = \wsclos{\sigma(\fA)} = \wotclos{\sigma(\fA)}$.  In general,
there is no analogue of this for operator algebras which are not
self-adjoint.  Indeed, it is possible to construct many examples of
pathology \cite{Wogen}.  On the other hand, the density theorem
is such a useful fact that it is worth seeking such a result whenever
possible.  In this section, we establish a density theorem for
regular representations of $\fA_n$.

Consider the following ``proof'' of the Kaplansky density theorem.
Consider the C*-algebra $\A$ sitting inside its double dual $\A^{**}$,
which is identified with the universal enveloping von Neumann algebra
$\W_u$ of $\A$.  
Any representation $\sigma$ of $\A$ extends uniquely to a normal
representation $\ol{\sigma}$ of $\W_u$ onto $\W = \sigma(\A)''$.  
Because this is a surjective $*$-homomorphism of C*-algebras, it is a
complete quotient map.  In particular, any element of the open ball of
$\W$ is the image of an element in the ball of $\W_u$.  Now by
Goldstine's Theorem, every element of the ball of $\A^{**}$ is the
weak-$*$ limit of a net in the ball of $\A$.  Mapping this down
into $\W$ by $\ol{\sigma}$ yields the result.

We call this a ``proof'' because the usual argument that $\W_u$ is
isometrically isomorphic to $\A^{**}$ requires the Kaplansky density
theorem.  Indeed, each state on $\A$ extends to vector state on
$\W_u$.  But the fact that all functionals on $\A$ have the same norm
on $\W_u$ follows from knowing that the unit ball is weak-$*$ dense in
the ball of $\W_u$.  It seems quite likely that the use of Kaplansky's
density theorem could be avoided, making this argument legitimate.

Nevertheless, we can use this argument to decide when such a result
holds in our context.  Moreover, in the C*-algebra context,
Kaplansky's theorem extends easily to matrices over the algebra
because they are also C*-algebras.  In our case, it follows from the
proof.

The double dual of $\fA_n$ may be regarded as a free semigroup
algebra, in the following way.  
We shall use it as a tool in the proof of the Kaplansky density theorem,
and we pause to highlight some of its features.

\begin{defn}
Regard $\fA_n$ as a subalgebra of $\E_n$.
Then the second dual $\fA_n^{**}$ is naturally identified with a weak-$*$
closed subalgebra of $\E_n^{**}$.
This will be called the \textit{universal free semigroup algebra}.
That this is a free semigroup algebra will follow from the discussion
below. We shall denote its structure projection by $P_u$.
\end{defn}

Denote by $j$ the natural inclusion of a Banach space into its double
dual. Then $j(\fA_n)$ generates $\E_n^{**}$ as a von Neumann algebra.  

If $\sigma$ is a $*$-representation of $\E_n$ on a Hilbert space $\H$,
then $\sigma$ has a unique extension to a normal $*$-representation
$\ol{\sigma}$ of $\E_n^{**}$ on the same Hilbert space $\H$.
Moreover, $\ol{\sigma}(\E_n^{**})$ is the von Neumann algebra
$\sigma(\E_n)''$ generated by $\sigma(\E_n)$.
  
Fix once and for all a universal representation $\pi_u$ of $\E_n$
acting on the Hilbert space $\H_u$ with the property that $\pi_u$
has infinite multiplicity, i.e. $\pi_u \simeq \pi_u^{(\infty)}$.
This is done to ensure that the \wot\ and weak-$*$ topologies coincide on
the universal von Neumann algebra $\fW_u = \pi_u(\E_n)''$.
Then $\ol{\pi_u}$ is a $*$-isomorphism of $\E_n^{**}$ onto $\fW_u$.
This carries $\fA_n^{**}$ onto the weak-$*$ closed subalgebra closure
$\fS_u$ of $\pi_u(\fA_n)$. This coincides with the \wot-closure, and thus
this is a free semigroup algebra.
Hence $\fA_n^{**}$ is a free semigroup algebra.

Since $\pi_u$ has infinite multiplicity and contains a copy of $\lambda$, 
its type L part is spanned by wandering vectors.  
So by Theorem~\ref{T:equiv_ac}, the range of $\ol{\pi_u}(P_u^\perp)$ is
$\vac(\pi_u)$.

\begin{prop}\label{universalAbsCont}
Let $\sigma$ be a representation of $\E_n$ and let $P_u\in\A_n^{**}$
be the the universal structure projection.  Then
$\ol{\sigma}(P_u^\perp)$ is the projection onto $\vac(\sigma)$.
\end{prop}

\Prf Consider the kernel of $\ol{\sigma}$.
There is a central projection $Q_\sigma\in \E_n^{**}$ such that
$\ker\ol{\sigma}= Q_\sigma\E_n^{**}$.  
Moreover, we may regard $\H$ as a closed subspace of $\H_u$ and
$\ol{\sigma}$ as given by multiplication by $Q_\sigma^\perp$, namely
$\ol{\sigma}(X) = Q_\sigma^\perp X |_{\H}$ for any $X\in\E_n^{**}.$

Let $M$ be the range of $\ol{\sigma}(P_u^\perp)$ and let $x\in
M$.  Then $x\in Q_\sigma^\perp P_u^\perp \H_u$, so $x$ belongs to
$\vac(\pi_u)$.  Thus for any $A\in\fA_n$, 
\[
 \ip{\sigma(A)x,x}=
 \ip{\pi_u(j(A)) Q_\sigma^\perp P_u^\perp x, Q_\sigma^\perp P_u^\perp} .
\]
As the range of $P_u^\perp$ consists of absolutely
continuous vectors, we see that this is an
absolutely continuous functional, so $x\in\vac(\sigma)$.

Conversely, if $x\in\vac(\sigma)$, then there exists an intertwiner
$X\in\X(\sigma)$ and $\zeta \in\ell^2(\bF_n^+)$ so that $x=X \zeta$. 
Observe that $Q_\sigma^\perp X$ belongs to $\X(\pi_u)$, hence $x\in
\vac(\pi_u)$.  Since the absolutely continuous part of $\pi_u$
coincides with the type L part of $\pi_u$, we conclude that $x\in
P_u^\perp\H_u\cap Q_\sigma^\perp \H_u$ and therefore
$\ol{\sigma}(P_u^\perp)x=x$.
\qed 

Since the type L part of a representation $\sigma$ is contained in the
absolutely continuous part, it follows that $\ol{\sigma}(P_u^\perp)
\geq P_\sigma^\perp.$
Notice that by the previous result, $\sigma$ is regular if and only if
$\ol{\sigma}(P_u^\perp) = P_\sigma^\perp$, where $P_\sigma$ is the
structure projection for $\fS_\sigma$.

\begin{prop}\label{regularrepn} Let $\sigma$ be a regular
$*$-representation of $\E_n$.  Then
$\wotclos{\sigma(\fA_n)}=\wsclos{\sigma(\fA_n)}$ and
$\ol{\sigma}(\fA_n^{**}) = \wsclos{\sigma(\fA_n)}.$ 
\end{prop}

\Prf Let $\fT:= \wsclos{\sigma(\fA_n)}$,
$\fS:=\wotclos{\sigma(\fA_n)}$ and let $\fW$ be the von-Neumann
algebra generated by $\sigma(\fA_n)$.  Let $P_\fT$ and $P_\fS$ be the
structure projections for $\wsclos{\sigma(\fA_n)}$ and
$\wotclos{\sigma(\fA_n)}$ respectively.  Then $P_\fT^\perp\geq
P_\fS^\perp$.  Since the absolutely continuous part of $\sigma$
contains the range of $P_\fT^\perp$, the regularity of $\sigma$ yields
that $P_\fT=P_\fS=\ol{\sigma}(P_u)$.  Hence $\fT = \fW P + \fT
P^\perp$ and $\fS = \fW P + \fS P^\perp$.  Moreover both $\fT P^\perp$
and $\fS P^\perp$ are canonically isomorphic to $\fL_n$ and the
isomorphisms agree on $\sigma(\fA_n)$.  Hence they are equal.
For typographical ease, write $P=P_\fT=P_\fS$.

Given $X\in\fS$, find $X'\in\E_n^{**}$ such that $\ol{\sigma}(X')=X$.
We may suppose that $X' = Q_\sigma^\perp X'$. This determines $X'$
uniquely, and $\ol{\sigma}$ is injective on $Q_\sigma^\perp \E_n^{**}$. 
By Proposition~\ref{universalAbsCont} and the regularity of $\sigma$,
$\ol{\sigma}(P_u) = P$.
So $\ol{\sigma}(P_u X' P_u^\perp) = P X P^\perp=0$, whence
$P_u X' P_u^\perp = 0$.
To see that $X'$ belongs to $\fA_n^{**}$, it remains to show that
$P_u^\perp X'P_u^\perp$ lies in $\fA_n^{**} P_u^\perp$, which is type L.
But $\fA_n^{**} P_u^\perp$ and $\fS P^\perp$ are both canonically
isometrically isomorphic to $\fL_n$, from which it is clear that
$\ol{\sigma}|_{\fA_n^{**}P_u^\perp}$ is an isomorphism onto
$\fS P^\perp$. 
\qed

We can now prove our Kaplansky-type theorem.

\begin{thm}\label{Kap}
Let $\sigma$ be a regular $*$-representation of $\E_n$.
Then the unit ball of $\sigma(\fA_n)$ is weak-$*$ dense in the unit ball
of $\wsclos{\sigma(\fA_n)}$, and the same holds for
$\fM_k(\sigma(\fA_n))$.

\end{thm}
\Prf Let $\fS:= \wsclos{\sigma(\fA_n)}=\wotclos{\sigma(\fA_n)}.$  We
first show that 
\begin{equation}\label{kersigma}
 \ker \ol{\sigma}|_{\fA_n^{**}}=\fA_n^{**}Q_\sigma P_u.
\end{equation}
To see this, notice that
$\ol{\sigma}|_{\fA_n^{**}P_u^\perp}$ is an isometric map of the type L
part of $\fA_n^{**}$ onto the type L part of $\fS$, that is,
$\ol{\sigma}$ maps $\fA_nP_u^\perp$ isometrically onto
$\fS P_{\sigma}^\perp$.
Therefore, if $X\in \fA_n^{**}$ and $\ol{\sigma}(X)=0$, then
$\ol{\sigma}(X)P_\sigma^\perp=0$, so that $XP_u^\perp=0$.  
As $X\in \ker\ol{\sigma}$, we find $X\in\fA_n^{**}Q_\sigma P_u$.  
The reverse inequality is clear, so \eqref{kersigma} holds.

Next we show that $\ol{\sigma}|_{\fA_n^{**}}$ is a complete quotient
map onto $\fS$.  For $X\in \fA_n^{**}$, we have 
\begin{align*}
\dist(X,\ker \ol{\sigma}|_{\fA_n^{**}}) &\leq \norm{X-XQ_\sigma P_u}\\
&=\norm{XQ_\sigma^\perp +XP_u^\perp Q_\sigma}\\
&=\max\{\norm{XQ_\sigma^\perp}, \norm{XP_u^\perp Q_\sigma}\}\\
&\leq \max\{\norm{XQ_\sigma^\perp}, \norm{XP_u^\perp}\}\\
&=\max\{\norm{\ol{\sigma}(X)}, \norm{\ol{\sigma}(X)P_\sigma^\perp}\}\\
&=\norm{\ol{\sigma}(X)}.
\end{align*}
The reverse inequality is clear, so that
$\norm{\ol{\sigma}(X)}=\dist(X,\ker\ol{\sigma}|_{\fA_n^{**}}).$ By
tensoring $Q_\sigma$ and $P_u$ with the identity operator on a
$k$-dimensional Hilbert space, the same argument holds for $X\in
M_k(\fA_n^{**})$ and the map $\sigma_k:=\sigma\otimes I_{\bC^k}$.
Thus $\ol{\sigma}|_{\fA_n^{**}}$ is a complete contraction.

Consider any element $T$ of the open unit ball of $\fS$.
Since the map of $\fA_n^{**}$ onto $\fS$ is a complete quotient map, 
there is a contraction $T_u \in \fA_n^{**}$ which maps onto $T$.
By Goldstine's Theorem, the unit ball of a Banach space is weak-$*$
dense in the ball of its double dual.
So select a net $A_\lambda$ in the ball of $\fA_n$ so that
$j(A_\lambda)$ converges weak-$*$ to $T_u$. 
Then evidently $\sigma(A_\lambda)$ converges weak-$*$ (and thus \wot) to $T$.
If one wants $\|A_\lambda\| \le \|T\|$, a routine modification will 
achieve this.

Because $\ol{\sigma}|_{\fA_n^{**}}$ is a complete contraction, the
same argument persists for matrices over the algebra as well.  \qed

\begin{lem}\label{KapWCnstIsTypeL}
If $\sigma$ is absolutely continuous and $\fS$ satisfies
Kaplansky's Theorem with a constant, then $\sigma$ is type L.
\end{lem}

\Prf 
As $\sigma$ is absolutely continuous, $\sigma \oplus \lambda$ is type L.
Let $\tau$ denote the weak-$*$ continuous homomorphism of $\fL_n$
into $\fS$ obtained from the isomorphism of $\fL_n$ with
$\fS_{\sigma \oplus \lambda}$ followed by the projection onto the
first summand.

Note that if $L$ is an isometry in $\fL_n$, then
$(\sigma \oplus \lambda)(L)$ is an isometry \cite[Theorem~4.1]{DKP}.
Hence $\tau(L)$ is an isometry as well.

Consider $\ker \tau$.
This is a weak-$*$ closed two-sided ideal in $\fL_n$.
If this ideal is non-zero, then the range of the ideal is spanned by
the ranges of isometries in the ideal \cite{DP2}.
In particular, the kernel would contain these isometries, contrary to
the previous paragraph.
Hence $\tau$ is injective.

Let $C$ be the constant in the density theorem for $\fS$.
If $T \in \fS$ and $\|T\| \le 1/C$, then there is a net
$A_i$ in the unit ball of $\fA_n$ such that $\sigma(A_i)$ converges
weak-$*$ to $T$. 
Drop to a subnet if necessary so that the net
$\lambda(A_i)$ converges weak-$*$ to an element $A \in \fL_n$.
Then $(\sigma \oplus \lambda)(A_i)$ converges weak-$*$ to $T\oplus A$.
Hence $\tau(A) = T$.
That means that $\tau$ is surjective, and hence is an isomorphism.

Now if $\sigma$ is not type L, then it is von Neumann type by
Proposition~\ref{P:exception}; and hence contains proper projections.   
But $\fL_n$ contains no proper idempotents \cite{DP1}; so this is
impossible. Therefore $\sigma$ must be type L.
\qed

\begin{thm}\label{KaplanskyIsRegular}
For a representation $\sigma$ of $\E_n$, the following
statements are equivalent. 
\begin{enumerate} 
\item The unit ball of
$\sigma(\fA_n)$ is \wot-dense in the ball of $\fS =
\wotclos{\sigma(\fA_n)}$.  i.e. Kaplansky's density theorem holds.
\item The \wot-closure of the unit ball of $\sigma(\fA_n)$ in $\fS =
\wotclos{\sigma(\fA_n)}$ has interior.  i.e. Kaplansky's density
theorem holds with a constant.  
\item $\sigma$ is regular.  
\end{enumerate}   
\end{thm}

\Prf 
(3) implies (1) follows from Theorem~\ref{Kap}.
That (1) implies (2) is obvious, so suppose (2) holds.  
If $\sigma$ is not regular, then it is von Neumann type by
Proposition~\ref{P:exception}; and $\sigma \simeq \sigma_a \oplus
\sigma_s$.  
Since Kaplansky holds with a constant, this persists for $\sigma_a$
because the \wot-closure does not change by dropping $\sigma_s$, it 
being the full von Neumann algebra already.
This contradicts Lemma~\ref{KapWCnstIsTypeL}.
\qed

\begin{defn}\label{D:sing}
A functional $\phi$ on $\fA_n$ is \textit{singular}
if it annihilates the type L part of $\fA_n^{**}$.
\end{defn}

\begin{prop}\label{P:sing}
For a functional $\phi$ on $\fA_n$ of norm $1$,
the following are equivalent:
\begin{enumerate}
\item $\phi$ is singular.
\item There is a regular representation $\sigma$ of $\E_n$
and vectors $x,y \in \H_\sigma$ with $x = P_\sigma x$
such that $\phi(A) = \ip{ \sigma(A) x, y }$.
\item $\dlim_{k\to\infty} \| \phi|_{ \fA_{n,0}^k} \| = 1$.
\end{enumerate}
If $\phi$ extends to a state on $\E_n\ ($i.e. $\phi(I)=1 )$, then $(3)$
is equivalent to \vspace{.5ex}\\ 
\strut\quad $(3')$\ $\| \phi|_{ \fA_{n,0}} \| = 1$.
\end{prop}

\Prf If $\phi\in\fA^*$, it is a weak-$*$ continuous functional on
$\fA_n^{**}$, so we may
represent it as a vector functional on $\H_u$, say
$\phi(A) = \ip{ \pi_u (A) x ,y }$.
Since $\phi$ annihilates the type L part, it does not change the
functional to replace $x$ by $P_u x$.
So (1) implies (2).

If (2) holds, then for every $A\in\fA_n$ we have 
$\phi(A)=\ip{\ol{\sigma}(j(A)P_u)x,y},$ which clearly annihilates the
type L part of $\fA^{**}$.  Thus (2) implies (1). 

If (1) holds, then $\phi(j(A)P_u^\perp)=0$, so $\phi(j(A))=\phi(j(A)P_u)$.
Now $\fA_n^{**}P_u = \bigcap_{k\geq 1} (\fA_{n,0}^{**})^k$, so that
$\| \phi|_{(\fA_{n,0}^{**})^k} \| =1$ for all $k \ge 1$.
It is easy to see that $(\fA_{n,0}^{**})^k = (\fA_{n,0}^k)^{**}$.
By basic functional analysis, a functional on a Banach space $X$ has the
same norm on the second dual.
Therefore $\| \phi|_{\fA_{n,0}^k} \| = 1$ for all $k \ge 1$.

If (3) holds, then there is a sequence $A_k$ in the ball of
$\fA_{n,0}^k$ so that $\dlim_{k\to\infty}\| \phi(A_k) \| = 1$.
Dropping to a weak-$*$ convergent subnet, we may assume that
this subnet converges to an element 
$A \in \bigcap_{k\ge1} (\fA_{n,0}^{**})^k = \fA_n^{**} P_u$.
Thus $\| \phi|_{\fA_n^{**} P_u} \| = 1$.

We claim that $\phi|_{\fA_n^{**} P_u^\perp}  = 0$.
If not, there is a norm one element $B = B P_u^\perp$ with
$\phi(B) = \ep > 0$.
Then 
\[
 \| A + \ep B \| = \| AA^* + \ep^2 BB^* \|^{1/2} \le
 (1+\ep^2)^{1/2} .
\]
But $\phi(A+\ep B) = 1 + \ep^2$, and thus $\| \phi \| > 1$.
This contradiction shows that $\phi$ annihilates the type L part,
and thus is singular.
So (3) implies (1).

Clearly (3) implies $(3')$.  Conversely, if $\phi$ extends to a state
on $\E_n$, it may be regarded as a normal state on $\E_n^{**}$ and
hence represented as $\phi(A) = \ip{ \sigma(A) \xi, \xi }$ where
$\sigma$ is obtained from the GNS construction. 
If $A\in \fA_{n,0}^{**}$  satisfies $1 = \|A\| = \phi(A)$, then
$\sigma(A)\xi = \xi$ is an eigenvalue.  
Therefore $\phi(A^k) = 1$ for all $k\ge1$, showing that 
$\| \phi|_{(\fA_{n,0}^{**})^k} \| =1$ for all $k \ge 1$.
Arguing as above establishes (3).
\qed

Here is a version of the Jordan decomposition.

\begin{prop}\label{decomp}
Every functional $\phi$ on $\fA_n$ splits uniquely as the sum of an
absolutely continuous functional $\phi_a$ and a singular one $\phi_s$.
Moreover 
\[
 \| \phi \| \le \| \phi_a \| + \| \phi_s \| \le \sqrt2 \| \phi \|
\]
and these inequalities are sharp.
\end{prop}

\Prf Set $\phi_a(A) = \phi( \pi_u(A) P_u^\perp)$ and 
$\phi_s = \phi( \pi_u(A) P_u)$.    
Clearly this is the desired decomposition.
For uniqueness, suppose that $\psi$ is both singular and absolutely
continuous.
Then $\|\psi\| = \dlim_{k\to\infty} \| \psi|_{\fA_{n,0}^k} \| = 0$.

Regard $\fA_n$ as a subalgebra of $\E_n$ and extend $\phi$ to a linear
functional (again called $\phi$) on $\E_n$ with the same norm.  Then
(using the GNS construction and the polar decomposition of functionals
on a C*-algebra) there exists a
$*$-representation $\sigma$ of $\E_n$ on a Hilbert space $\H_\sigma$
and vectors $x,y\in \H_\sigma$ with $\|x\|\,\|y\| =\norm{\phi}$ so
that for every $A\in\E_n$, $\phi(A) = \ip{ \sigma(A) x,y }.$ Therefore
for $A \in\fA_n$, we have $\phi_a(A) = \ip{ \sigma(A)
\ol{\sigma}(P_u^\perp) x , y }$ and $\phi_s(A) = \ip{ \sigma(A)
\ol{\sigma}(P_u) x , y }$.  Hence
\begin{align*}
 \| \phi_a \| + \| \phi_s \| &\le \| \ol{\sigma}(P^\perp) x \| \, \| y
 \| + \| \ol{\sigma}(P_u) x \| \, \| y \| \\&\le \sqrt2 \big( \|
 \ol{\sigma}(P^\perp) x \|^2 + \| \ol{\sigma}(P) x \|^2 \big)^{1/2}
 \|y\| = \sqrt2 \| \phi \| .
\end{align*}
The example following will show that the $\sqrt2$ is sharp.
\qed

\begin{eg}
Consider the atomic representation $\sigma_{1,1}$ on 
$\bC \xi_* \oplus \Hn$ given by $S_1 \xi_* = \xi_*$ and
$S_2 \xi_* = \xi_\mt$; and $S_i|_{\Hn} = L_i$.
Set 
\[
 \phi(A) = 
 \ip{ \sigma_{1,1}(A) (\xi_* + \xi_\mt)/\sqrt2, \xi_\mt } .
\]
Then $\fS_\sigma$ contains $A = \xi_\mt \xi_*^* /\sqrt2 +
(I-\xi_* \xi_*^*) /\sqrt2$ and $\phi(A)=1$.
So we see that $\|\phi\|=1$.

On the other hand, 
\[
 \phi_s(A) = \ip{ \sigma_{1,1}(A) \xi_*/\sqrt2, \xi_\mt } \qand 
 \phi_a(A) = \ip{ \sigma_{1,1}(A) \xi_\mt/\sqrt2, \xi_\mt }
\]
both have norm $1/\sqrt2$.  
So $\|\phi_s\| + \|\phi_a\| = \sqrt2$.
\end{eg}

\begin{ques}
Let $S$ be the unilateral shift and consider the representation of
$\fA_2$ obtained from the minimal isometric dilation of $A_1 =
S/\sqrt2$ and $A_2 = (S+P_0)^*/\sqrt2$.  The weak-$*$ closed
self-adjoint algebra generated by $A_1$ and $A_2$ is all of $\B(\H)$.
Therefore this representation is either dilation type with $P \fS P =
\B(\H)$ or it is type L, depending on whether the functional $\phi =
e_0 e_0^*$ is singular or absolutely continuous.  To check. it
suffices to determine whether $\phi$ has norm 1 or less on
$\fA_{n,0}$.  We would like to know which it is.
\end{ques}

\begin{ques} Charles Read has given an example of a representation of
$\fA_2$ such that $\B(\H) = \wsclos{\sigma(\fA_2)}$.  
Is Read's example singular or absolutely continuous?
Again it suffices to take any convenient state on $\B(\H)$
and estimate its norm on $\fA_{2,0}$ as equal to 1 or strictly less. 
\end{ques}

We provide an example of how the density theorem can be used to
establish an interpolation result for finitely correlated
presentations.  Such representations are obtained from a row
contraction of matrices $A =
\begin{bmatrix}A_1&\dots&A_n\end{bmatrix}\in M_{1,n}(M_k(\bC))$ by
taking the minimal isometric dilation \cite{Fr1,Bun,Pop_diln}.  These
representations were classified in \cite{DKS}.  The structure
projection $P$ has range equal to the span of all $\{A_i^*\}$
invariant subspaces on which $A$ is isometric.  In particular, it is
finite rank.  Also, the type L part is a finite multiple, say $\alpha$,
of the left regular representation.  Thus elements of the free
semigroup algebra $\fS$ have the form
$\begin{bmatrix}X&0\\Y&Z^{(\alpha)}\end{bmatrix}$ where $X$ and $Y$
lie in $P\fW P$ and $P^\perp \fW P$ respectively, and $Z \in \fL_n$,
where $\fW$ is the von Neumann algebra generated by $\fS$.

\begin{thm}\label{interpolate}
Let $\sigma$ be a finitely correlated representation.
If $A \in \fS_\sigma$ has $\|A\|<1$ and $k \in \bN$,
then there is an operator $B \in \fA_n$ so that $\sigma(B)P=AP$ and the Fourier
series of $B$ up to level $k$ agree with the coefficients of $AP^\perp$.
\end{thm}

\Prf Fix $\ep < 1 - \|A\|$.
Let $Q_k\in \B(\Hn)$ denote the projection onto $\spn\{ \xi_w : |w|
\le k\}$.  

Identify $\fS P^\perp$ with $\fL_n^{(\alpha)}$ and find $C\in \fL_n$ so
that  $AP^\perp =C^{(\alpha)}$.

Since the weak and strong operator topologies have the same closed
convex sets, the density theorem implies that there exists a sequence
$\{L_k\}$ in $\fL_n$ so that $\norm{L_k}< 1-\eps$ and $A=\sot\lim
\sigma(L_k)$.  Recalling that $P$ and $Q_k$ are finite rank, we
conclude that there exists $B_1\in \fA_n$ so that
\[
 \| ( A - \sigma(B_1) ) P \| + \| Q_k (C - B_1) \| < \ep/2 .
\]
By \cite[Corollary~3.7]{DP3}, there is an element $C_1 \in \fL_n$ so
that $Q_k C_1 = Q_k (C - B_1)$ and $\|C_1\| = \| Q_k (C - B_1) \|$.
Hence the element of $\fS$ defined by 
$A_1 = ( A - \sigma(B_1) ) P + C_1^{(\alpha)} P^\perp$
satisfies $\| A_1 \| < \ep/2$.

Now choose $B_2 \in \fA_n$ so that $\|B_2\| < \ep/2$ and
\[
 \| ( A_1 - \sigma(B_2) ) P \| + \| Q_k (C_1 - B_2) \| < \ep/4 .
\]
Proceed as above to define $C_2 \in \fL_n$ so that
$Q_k C_2 = Q_k (C_1 - B_2)$ and $\|C_2\| = \| Q_k (C_1 - B_2) \|$;
and then define $A_2 = ( A_1 - \sigma(B_2) ) P + C_2^{(\alpha)} P^\perp$
satisfying $\| A_2 \| < \ep/4$.

Proceeding recursively, we define $B_j$ for $j \ge 1$ so that 
$B = \sum_{j\ge1} B_k$ is the desired approximant.
\qed

%%%%%%    %%%%%%
\section{Constructive examples of Kaplansky}

In this section, we give a couple of examples where we were able to
construct the approximating sequences more explicitly.
We concentrate on exhibiting the structure projection $P$ as a limit of
contractions.
It is then easy to see that the whole left ideal $\fW P$ has the same
property by applying the C*-algebra Kaplansky theorem.
We do not have an easy argument to show that one can extend this to the
type L part without increasing the constant. 

\begin{prop}\label{oldKap} 
Let $\fS$ be the free semigroup algebra generated by isometries
$S_1,\dots, S_n$; and let $\fA$ be the norm closed algebra that they
generate.  Let $P\in\fS$ be the projection given by the Structure
Theorem. If $P\ne I$ is the \wot-limit of a sequence in $\fA$ of norm at
most $r$, then the \sot-closure of the $r$-ball of $\fA_0^k$
contains $\fS P$ for all $k \ge0$.
\end{prop}

\Prf Since $P\neq I$, $\fS$ has a type L part.  
Let $\Phi$ be the canonical surjection of $\fS$ onto $\fL_n$ with
$\Phi(S_i)=L_i$ \cite[Theorem~1.1]{DKP}.   
Recall that the kernel of $\Phi$ is $\bigcap_{k=1}^\infty \fS_0^{k} = \fW
P$. Since the weak and strong operator topologies have the same closed
convex sets, we may suppose that the sequence in $\fA$ converges
to $P$ strongly. 
In particular, the restriction of this sequence to the type L part
converges strongly to $0$.
Hence the Fourier coefficients each converge to 0.
Thus a minor modification yields a sequence $A_k \in \fA_0^k$ of norm at
most $r$ converging \sot\ to $P$.

If $T$ lies in the unit ball of $\fW P$, then by the usual Kaplansky
density theorem, there is a sequence $B_k$ in the unit ball of $\ca(S)$
which converges \sot\ to $T$.
We may assume that $B_k$ are polynomials in $S_i, S_i^*$ for $1 \le i \le n$
of total degree at most $k$.
Then observe that $B_k A_{2k}$ lies in $\fA_{n,0}^k$, and converges \sot\
to $TP = T$.
\qed

Our first example is a special class of finitely correlated representations
which are obtained from dilating multiples of unitary matrices.

\begin{thm}\label{unitary_Kap}
Suppose that $U_i$ for $1 \le i \le n$ are unitary matrices in
$\B(\V)$, where $\V$ has finite dimension $d$, and that $\alpha_i$
are non-zero scalars so that  $\sum_{i=1}^n |\alpha_i|^2 = 1$. 
Let $S_i$ be the joint isometric dilation of $A_i = \alpha_i U_i$
to a Hilbert space $\H$.  Let $\fS$ be the free semigroup algebra
that they generate; and let $\fA$ denote the norm-closed algebra. 
Then the projection $P = P_\V$ is the projection that occurs in
the Structure Theorem, and there is a sequence of contractions in
$\fA$ which converges \sot\ to $P$.
\end{thm}

\begin{lem}\label{unitary_algebra}
If $\U$ is a set of unitary matrices in $\fM_d$, then the closure
of the set of all non-empty words in elements of $\U$ is a
subgroup of the unitary group $\fU$ and the algebra generated by
$\U$ is a C*-algebra.
\end{lem}

\Prf The closure $\G$ of words in $\U$ is multiplicative and
compact. Any unitary matrix $U$ is diagonalizable with finite
spectrum. A routine pigeonhole argument shows that there is a
sequence $U^{n_i}$ which converges to $I$, and thus $U^{n_i-1}$
converges to $U^{-1}$.  It follows that $\G$ is a group.
It is immediate that the algebra generated by $\U$ contains $\U^*$
and thus is self-adjoint.
\qed

\noindent\textbf{Proof of Theorem~\ref{unitary_Kap}.\ }
{From} the Lemma, we see that the algebra generated by $\{ A_i^* \}$
is self-adjoint, and thus the space $\V$ is the span of its
minimal $A_i^*$ invariant subspaces.  From the Structure Theorem
for finitely correlated representations \cite{DKS}, we deduce that
$P=P_\V$ is a projection in $\fS$ and that 
$\fS = \fW P + \fS P^\perp$, where $\fW$ is the von Neumann algebra
generated by $\fS$, $P^\perp \H$ is invariant, and  $\fS P^\perp$
is a (finite) ampliation of $\fL_n$.

Consider the space $\fX$ consisting of all infinite words $x = i_1 i_2 i_3
\dots$ where $1 \le i_j \le n$ for $j \ge1$.  
This is a Cantor set in the product topology.
Put the product measure $\mu$ on $\fX$ obtained from the measure
on $\{1,\dots,n\}$ which assigns mass $|\alpha_i|^2$ to $i$.

Fix $\ep>0$.  Since the closed semigroup $\G$ generated by
$\{U_i\}$ is a compact group by Lemma~\ref{unitary_algebra}, one
may choose a finite set $\S$ of non-empty words which form an
$\ep$-net (in the operator norm).  
Let $N$ denote the maximum length of these words.
Then we have the following consequence: given any word 
$w = i_1 \dots i_k$, there is a word $v = j_1 \dots j_l$ in $\S$ with 
$l \le N$ so that $U_{wv} = U_{i_1} \dots U_{i_k} U_{j_1} \dots U_{j_l}$
satisfies $\| U_{wv} - I \| < \ep$.

Recursively determine a set $\W$ of words so that $S_w$ have
pairwise orthogonal ranges and $\| U_w - I \| < \ep$ for $w \in
\W$ as follows:  start at an arbitrary level $k_0$ and take all
words $w$ with $|w| = k_0$ such that $\| U_w - I \| < \ep$.
If a set of words of length at most $k$ has been selected, add to
$\W$ those words of length $k+1$ which have ranges orthogonal to
those already selected and satisfy $\| U_w - I \| < \ep$.

We claim that $\sotsum_{w \in \W} S_w S_w^* = I$.
The argument is probabilistic.  
Let $\delta = \min \{ |\alpha_i|^2 \}$.
Associate to $w$ the subset $X_w$ of all infinite words in $\fX$
with $w$ as an initial segment.
By construction, the sets $X_w$ are pairwise disjoint clopen
sets for $w \in \W$ with measure $|\alpha_w|^2$, where we set
$\alpha_w = \prod_{t=1}^k \alpha_{i_t}$.
Verifying our claim is equivalent to showing that
$\bigcup_{w\in\W} X_w$ has measure 1.
Consider the complement $Y_k$ of $\bigcup_{w\in\W,\ |w|\le k} X_w$.
This is the union of certain sets $X_w$ for words $w$ of length $k$.
For each such word, there is a word $v\in \S$ so that 
$\| U_{wv} - I \| < \ep$.  Now $S_{wv}$ has range contained in the range
of $S_w$, which is orthogonal to the ranges of words in $\W$ up to level
$k$.
It follows from the construction of $\W$ that there will be a word 
$w' \in \W$ so that $w'$ divides $wv$.
As a consequence, $Y_{k+N}$ has measure smaller than $Y_k$ by a
factor of at most $1- \delta^N$ because for each interval $X_w$
in $Y_k$, there is an interval $X_{wv}$ which is in the complement,
and its measure is at least $\delta^N \mu(X_w)$.
Therefore $\dlim_{k\to\infty} \mu(Y_k) = 0$.

Choose a finite set $\W_0 \subset \W$ so that 
\[
 r := \mu( \bigcup_{w\in\W_0} X_w ) > 1 - \ep .
\]
Define $T = \sum_{w\in\W_0} \ol{\alpha_w} S_w$ in $\fA$.
Note that 
\[
 \|T\|^2 = \sum_{w\in\W_0} |\alpha_w|^2 = r < 1 .
\]
Observe that $P S_w = \alpha_w U_w$.
Define a state $\tau$ on $\B(\H)$ as the normalized trace of the
compression to $\V$.
Since $\|U_w - I \| < \ep$, it follows that 
$| \tau(U_w) - 1 | < \ep$.  Compute 
\begin{align*}
  \re \tau (T) &= \sum_{w\in\W_0} \ol{\alpha_w} \re \tau(S_w) 
   = \sum_{w\in\W_0} |\alpha_w|^2 \re \tau(U_w) \\
  &\ge \sum_{w\in\W_0} |\alpha_w|^2 (1-\ep) = r ( 1 - \ep)
  > ( 1 - \ep)^2 .
\end{align*}

By taking $\ep = 1/k$ and $k_0=k$ in the construction above, we
obtain a sequence $T_k$ of polynomials $T_k \in \fA_0^k$ which
are contractions and $\dlim_{k\to\infty} \tau(T_k) = 1$.
It follows that there is a subsequence which converges \wot\ to a
limit $T \in \fS$ which lies in $\bigcap_{k\ge1} \fS_0^k = \fW P$.
Moreover $\|T\| \le 1$ and $\tau(T) = 1$.
The only contraction in $\B(\V)$ with trace 1 if the identity, and
therefore the compression $PT = P$.  As $T$ is contractive, we
deduce that $P^\perp T = P^\perp TP = 0$, whence $T=P$.
As the \sot\ and \wot-closures of the balls are the same, there is
a sequence in the convex hull of the $T_k$'s which converges to
$P$ strongly.
\qed

Our second constructive example is the set of atomic representations
introduced in \cite{DP1}.
To analyze these, we will need some of Voiculescu's theory of free
probability.

\begin{thm}\label{Kap_atomic}
If $\fS$ is an atomic free semigroup algebra, then the structure
projection is a \sot-limit of contractive polynomials in the
generators.  
\end{thm}

It is convenient for our calculation to deal with certain norm
estimates in the free group von Neumann algebra.  We thank Andu
Nica for showing us how to handle this free probability machinery.

\begin{lem}\label{Nica} 
Let $p$ and $q$ be free proper projections of trace $\alpha \le 1/2$
in a finite von Neumann algebra $(\fM, \tau)$. 
Then  $\|pqp\| = 4 \alpha (1-\alpha)$.
\end{lem}

\Prf
Given $a \in \fM$, form the power series $M_a(z) =
\dsum_{n\ge1} \tau(a^n) z^n$.  
In particular, $M_p(z) = M_q(z) = \alpha z (1-z)^{-1}$.\vspace{.5ex}
Voiculescu's S-transform is given by
$S_a(\zeta) = (1+\zeta) \zeta^{-1} M_a^{<-1>}(\zeta)$
where $M_a^{<-1>}$ denotes the inverse of $M_a$ under composition.
So 
\[
 S_p(\zeta) = S_q(\zeta) = \dfrac{1+\zeta}{\zeta}
\dfrac{\zeta}{\alpha + \zeta} = \dfrac{1+\zeta}{\alpha + \zeta} .
\]
By \cite[Theorem 2.6]{V}, the
S-transform is multiplicative on free pairs. Hence 
\[
 S_{pq} = \dfrac{(1+\zeta)^2}{(\alpha + \zeta)^2}
        = \dfrac{1+\zeta}{\zeta} M_{pq}^{<-1>}(\zeta) .
\]
Observe that $\tau((pqp)^n) = \tau((pq)^n)$ and so $M_{pq} =
M_{pqp}$.  Thus 
\[
 M_{pqp}^{<-1>}(\zeta) = 
 \dfrac{\zeta(1+\zeta)}{(\alpha + \zeta)^2} .
\]

We obtain the quadratic equation
\[
 (1-z) M_{pqp}^2(z) + (1 - 2 \alpha z) M_{pqp}(z) - \alpha^2 z = 0 .
\]
Solving, we obtain
\[
 M_{pqp}(z) =
 \frac{2 \alpha z - 1 + \sqrt{1 - 4\alpha (1- \alpha)z} }{2(1-z)}
\]
Since $M_{pqp}(0)=0$, we must choose an appropriate branch of the
function $(1 - 4\alpha (1-\alpha)z)^{1/2}$.  
This may be defined on the complement of the line segment 
\[ \{ z : \re z = 1/4\alpha (1- \alpha) \AND \im z \le 0 \} \]
and takes positive real values on real numbers
$x < 1/4\alpha (1- \alpha)$.
An easy calculation shows that the singularity at $z=1$ is
removable.

The power series for $M_{pqp}$ converges on the largest disk on
which it is analytic.  The branch point occurring at 
$z = 1/4\alpha (1- \alpha)$ is the only obstruction, and thus
the radius of convergence is $1/4\alpha (1- \alpha)$.
On the other hand, from Hadamard's formula, the reciprocal of the
radius of convergence is 
\[
 4\alpha (1- \alpha) = \limsup_{k\to\infty} \tau((pqp)^k)^{1/k} = 
 \| pqp \| .
\]
\upqed

\begin{cor}\label{free_sum}
Let $U_i$ for $1 \le i \le n$ denote the generators of the free
group von Neumann algebra, and let $P_i$ be spectral projections
for $U_i$ for sets of measure at most $\alpha \le 1/2$.  
Then $\| \sum_{i=1}^n P_i \| \le 1 + 2n^2 \sqrt{\alpha}$.
\end{cor}

\Prf By Lemma~\ref{Nica}, we have $\|P_i P_j\| = 
\|P_i P_j P_i\|^{1/2} \le 2 \sqrt{\alpha}$ for $i \ne j$.
If $\| \sum_{i=1}^n P_i \| = 1 + x$, then
\[
 (1+x)^2 = \| \sum_{i=1}^n P_i + \sum_{i\ne j} P_i P_j \|
 \le 1+x + n(n-1) 2 \sqrt{\alpha}.
\]
Hence $x \le 2n^2 \sqrt{\alpha}$ as claimed.
\qed

Recall from \cite{DP1} the atomic representation $\sigma_{u,\lambda}$
determined by a primitive word $u = i_1\dots i_d$ in $\Fn$ and a scalar
$\lambda$ in $\bT$.  
Define a Hilbert space $\H_u \simeq \bC^d \oplus \Hn^{d(n-1)}$
with orthonormal basis $\zeta_1,\dots,\zeta_d$ for $\bC^d$ and index
the  copies of $\Hn$ by $(s,j)$, where 
$1 \le s \le d$, $1 \le j \le n$ and $j \ne i_s$,
with basis $\{\xi_{s,j,w}: w \in \Fn \}$. 
Define a representation $\sigma_{u,\lambda}$ of $\Fn$ and
isometries $S_i = \sigma_{u,\lambda}(i)$ by
\begin{alignat*}{2}
S_i \zeta_s &= \zeta_{s-1} &&\qif i=i_s,\,s>1\\ 
S_i \zeta_1 &= \lambda \zeta_d &&\qif i=i_1\\ 
S_i \zeta_s &= \xi_{s,i,\mt} &&\qif i\ne i_s\\
S_i \xi_{s,j,w} &= \xi_{s,j,iw} &&\qforal i,s,j,w
\end{alignat*}
For our purposes, we need to observe that the vectors
$\zeta_1,\dots,\zeta_d$ form a ring which is cyclically permuted by the
appropriate generators $S_{i_s}$; and all other basis vectors are wandering.
The projection $P$ in the structure theorem is the projection onto $\bC^d$.

\begin{lem}\label{atomic1}
Let $u$ be a primitive word and let $\lambda \in \bT$.  
Let $\fS$ be the atomic free semigroup algebra corresponding to the
representation $\sigma_{u,\lambda}$.  
Then the projection $P$ from the Structure Theorem is the limit of
contractive polynomials in the generators.
\end{lem}

\Prf   
Let $u_s$ denote the cyclic permutations of $u$ for $1 \le s \le d$
satisfying $S_{u_s} \zeta_s = \lambda \zeta_s$.
As in the proof \cite[Lemma~3.7]{DP1} of the classification of atomic
representations, there is a sequence of the form  $A_{k,s} = p_k(S_{u_s})$
which converge \sot\ to the projections $\zeta_s \zeta_s^*$ where    
$p_k(x) = x^k q_k(x)$ are polynomials with
$\|p_k\|_\infty = 1 = p_k(\lambda)$.  

It is routine to choose such polynomials $p_k$ with the added
stipulation that there is an open set $\V_k$ of measure $1/k$
containing the point $\lambda$ so that 
$\| p_k \upchi_{ \bT \bsl \V_k } \| < 1/k$.
We now consider the elements $A_k = \sum_{s=1}^d A_{k,s}$.
Clearly the sequence $A_k$ converges \sot\ to the projection $P$.
So it suffices to establish that $\dlim_{k\to\infty} \|A_k\| = 1$.

As $\fA$ has a unique operator algebra structure independent of the
representation of $\O_n$, polynomials in the isometries $S_{u_s}$
may be replaced by the corresponding polynomials in the generators
$L_s$ of the left regular representation for the free semigroup
$\bF_d^+$.
The isometries $L_s$ have pairwise orthogonal ranges for
distinct $s$, and consequently the operators $A_{k,s}$ have
orthogonal ranges.  
Hence the norm of $A_k$ equals the norm of the column operator with
entries $A_{k,s}$. 
Now it is evident that the left regular representation of
$\bF_d^+$ may be obtained as the restriction of the left regular
representation of the free group $\bF_d$ to an invariant subspace.
So the norm is increased if $L_s$ are replaced by the generators
$U_s$ of the free group von Neumann algebra.
Thus the norm of $A_k$ is dominated by the column vector with
entries $p_k(U_s)$.

Let $Q_s$ denote the spectral projection of $U_s$ for the open set
$\V_k$.  
Then 
\begin{align*}
  \begin{bmatrix}p_k(U_1)\\ \vdots\\ p_k(U_d) \end{bmatrix} &=
  \diag\{p_k(U_1), \dots, p_k(U_d) \}
  \begin{bmatrix}Q_1\\ \vdots\\ Q_d \end{bmatrix} + 
  \begin{bmatrix}p_k(U_1)Q_1^\perp\\ \vdots\\ p_k(U_d)Q_d^\perp \end{bmatrix}.
\end{align*}
Since $\|p_k(U_s)\| = 1$, the norm of the first term is at most
\[
 \| Q_1 + \dots + Q_d \|^{1/2} < 1 + 2 d^2/ \sqrt k
\]
by Corollary~\ref{free_sum}.
The second term is dominated by 
\[
 \sqrt d\,\, \| p_k \upchi_{ \bT \bsl \V_k } \|_\infty < \sqrt d / k .
\]
Hence the norms converge to 1 as claimed.
\qed

It is now only a technical exercise to show how one may use similar
arguments to combine sequences corresponding to finitely many points on the
circle for a given word $u$, and to deal with finitely many such words at
once.  The inclusion of summands of type L such as the atomic
representations of inductive type does not affect things since these
sequences are already converging strongly to 0 on the wandering subspaces
of these atomic representations.  Details are omitted.

\end{document}